\documentclass[aos,seceqn,nameyear,dvips]{arximspdf}
\usepackage{dcolumn}
\usepackage{graphicx}

\doi{10.1214/10-AOS848}
\volume{39}
\issue{1}
\pubyear{2011}
\firstpage{613}
\lastpage{642}

\makeatletter

\newcolumntype{d}[1]{D{.}{.}{#1}}

\newtheorem{theorem}{Theorem}[section]
\newtheorem{lemma}{Lemma}

\newproclaim{remark}{Remark}

\newtheorem{corollary}{Corollary}

\newproclaim{definition}{Definition}

\newproclaim{step}{Step}
\newproclaim{Condition}{Condition}

\newcommand{\FDP}{\operatorname{FDP}}
\newcommand{\FDR}{\operatorname{FDR}}
\newcommand{\FWER}{\operatorname{FWER}}
\newcommand{\Betadist}{\operatorname{Beta}}
\newcommand{\RR}{\mathbb{R}}
\newcommand{\ZZ}{\mathbb{Z}}
\newcommand{\ID}{\mathrm{I}}
\newcommand{\versus}{{\mbox{versus}}}
\newcommand{\calV}{\mathcal{V}}
\newcommand{\chiK}{\mathbb{K}}
\newcommand{\bc}{\mathrm{bc}}
\newcommand{\AFNI}{\operatorname{AFNI}}
\newcommand{\FSL}{\operatorname{FSL}}
\newcommand{\LIP}{\operatorname{LIP}}
\newcommand{\bh}{{\mathbf{h}}}
\newcommand{\bS}{\mathbf{S}}
\newcommand{\br}{{\mathbf{r}}}
\newcommand{\fMRI}{\operatorname{fMRI}}

\makeatother

\begin{document}
\begin{frontmatter}

\title{Multiple testing via $\mbox{FDR}_{L}$ for large-scale
imaging~data}
\runtitle{$\mbox{FDR}_L$ for large-scale imaging data}

\begin{aug}
\author[A]{\fnms{Chunming} \snm{Zhang}\corref{}\thanksref{t1}\ead[label=e1]{cmzhang@stat.wisc.edu}},
\author[B]{\fnms{Jianqing} \snm{Fan}\thanksref{t2}\ead[label=e2]{jqfan@princeton.edu}} and
\author[C]{\fnms{Tao} \snm{Yu}\ead[label=e3]{stayt@nus.edu.sg}}
\runauthor{C. Zhang, J. Fan and T. Yu}
\affiliation{University of Wisconsin-Madison, Princeton
University\break
 and National University of Singapore}
\address[A]{C. Zhang\\
Department of Statistics\\
University of Wisconsin\\
Madison, Wisconsin 53706\\
USA\\
\printead{e1}}
\address[B]{J. Fan\\
Department of Operation Research\\
\quad and Financial Engineering \\
Princeton University \\
Princeton, New Jersey 08544\\
USA \\
\printead{e2}}
\address[C]{T. Yu\\
Department of Statistics\\
\quad and Applied Probability \\
National University of Singapore \\
Singapore 117546 \\
\printead{e3}}
\end{aug}

\thankstext{t1}{Supported by NSF Grant DMS-07-05209 and Wisconsin Alumni
Research Foundation.}
\thankstext{t2}{Supported by NIH R01-GM072611 and NSF Grant DMS-07-14554.}

\received{\smonth{5} \syear{2010}}
\revised{\smonth{8} \syear{2010}}

%
\begin{abstract}
The multiple testing procedure plays an important role in detecting the
presence of spatial signals for large-scale imaging data. Typically,
the spatial signals are sparse but clustered. This paper provides
empirical evidence that for a range of commonly used control levels,
the conventional $\FDR$ procedure can lack the ability to detect
statistical significance, even if the $p$-values under the true null
hypotheses are independent and uniformly distributed; more generally,
ignoring the neighboring information of spatially structured data will
tend to diminish the detection effectiveness of the $\FDR$ procedure.
This paper first introduces a scalar quantity to characterize the
extent to which the ``\textit{lack of identification phenomenon}''
($\LIP$) of the $\FDR$ procedure occurs. Second, we propose a new
multiple comparison procedure, called $\FDR_L$, to accommodate the
spatial information of neighboring $p$-values, via a local aggregation
of $p$-values. Theoretical properties of the $\FDR_L$ procedure are
investigated under weak dependence of $p$-values. It is shown that the
$\FDR_L$ procedure alleviates the $\LIP$ of the $\FDR$ procedure, thus
substantially facilitating the selection of more stringent control
levels. Simulation evaluations indicate that the $\FDR_L$ procedure
improves the detection sensitivity of the $\FDR$ procedure with little
loss in detection specificity. The computational simplicity and
detection effectiveness of the $\FDR_L$ procedure are illustrated
through a real brain fMRI dataset.
\end{abstract}

%
\begin{keyword}[class=AMS]
\kwd[Primary ]{62H35}
\kwd{62G10}
\kwd[; secondary ]{62P10}
\kwd{62E20}.
\end{keyword}
\begin{keyword}
\kwd{Brain fMRI}
\kwd{false discovery rate}
\kwd{median filtering}
\kwd{$p$-value}
\kwd{sensitivity}
\kwd{specificity}.
\end{keyword}

\end{frontmatter}

\section{Introduction} \label{sec-1}

In many important applications, such as astrophysics, satellite
measurement and brain imaging,
the data are collected at spatial grid points, and a large-scale
multiple testing procedure
is needed for detecting the presence of spatial signals.
For example, functional magnetic resonance imaging (fMRI) is a recent
and exciting imaging
technique that allows
investigators to determine which areas of the brain are involved in a
cognitive task.
Since an fMRI dataset contains time-course measurements over voxels,
the number of which is typically of the order of $10^4$--$10^5$, a
multiple testing procedure plays an important role in
detecting the regions of activation. Another example of important
application of multiple testing
is to the diffusion tensor imaging, which intends to identify brain
white matter regions
[\citet{Bihan01}].

In the seminal work, \citet{Wetal02} proposed a Gaussian random field
method which approximates the family-wise error rate ($\FWER$) by
modeling test
statistics
over the entire
brain as a Gaussian random field. It has been found to be
conservative in some cases [\citet{NH03}]. \citet{NH03} also discussed the use of permutation tests and their
simulation studies showed that permutation tests tended to be more
sensitive in
finding activated regions.
The false discovery rate ($\FDR$) approach has become increasingly popular.
The conventional $\FDR$ procedure offers
the advantage of
overcoming the conservativeness drawback of $\FWER$,
requiring fewer assumptions than random field based methods and
being computationally less intensive than permutation tests.

Nevertheless, in practical applications to imaging data with a spatial
structure,
even if the $p$-values corresponding to the true null hypotheses are
independent and uniformly distributed,
the conventional $\FDR$ procedure may lack the ability to detect
statistical significance,
for a range of commonly used control levels $\alpha$.
It will be seen, in the left panels of Figure \ref{Figure-2}, that
the $\FDR$ procedure for a 2D simulated data declares only a couple of
locations to be significant for $\alpha$
ranging from $0$ to about $0.4$.
That is, even if we allow $\FDR$ to be controlled at the level $40\%$,
one cannot reasonably well identify
significant sites.
The empirical evidence provided above for the standard $\FDR$
procedure is not pathological. Indeed, similar phenomena arise from commonly
used signals plus noise models for imaging data,
as will be exemplified by
extensive studies in Section \ref{sec-4.2}.
In statistical literature, while some useful finite-sample and
asymptotic results [\citet{STS04}]
have been established for the $\FDR$ procedure, the results could not
directly quantify
the loss of power and ``\textit{lack of identification phenomenon}''
($\LIP$).

More generally, for spatially structured imaging data, the significant
locations are typically sparse,
but clustered rather than scattered. It is thus anticipated that a
location and its adjacent neighbors
fall in a similar type of region, either significant (active) or
nonsignificant (inactive).
As will be seen in the simulation studies
(where the $\LIP$ does not occur)
of Section \ref{sec-5},
the existing $\FDR$ procedure
tends to be less effective in detecting significance.
This lack of detection efficiency is due to the information of
$p$-values from adjacent neighbors not having been
fully taken into account.
Due to the popularity of the $\FDR$ procedure in research practices,
it is highly desirable
to embed the spatial information of imaging data into the $\FDR$ procedure.

This paper aims
to quantify the $\LIP$
and
to propose a new multiple testing procedure, called $\FDR_L$, for
imaging data, to
accommodate the spatial information of neighboring $p$-values, via a
local aggregation of $p$-values.
Main results are given in three parts.
\begin{itemize}
\item In the first part, statistical inference for the null
distribution of locally aggregated $p$-values is studied.
See Method I proposed in Section \ref{sec-3.2}
and Method II in Section \ref{sec-3.3}.

\item In the second part, asymptotic properties of the $\FDR_L$
procedure are investigated under
weak dependence (to be defined in Section \ref{sec-4.1}) of
$p$-values. See Theorems \ref{Thm-1}--\ref{Thm-3}.

\item The third part intends to
provide a more in-depth discussion of why the $\LIP$ occurs and the
extent to which the $\FDR_L$ procedure alleviates
the $\LIP$.
In particular, we introduce a scalar $\alpha_\infty$ to quantify the
$\LIP$:
the smaller the $\alpha_\infty$,
the smaller control level can be adopted without encountering $\LIP$;
$\alpha_\infty=0$ rules out the possibility of the $\LIP$.
In the particular case of i.i.d. $p$-values, Theorem \ref{Thm-4}
provides verifiable
conditions under which $\alpha_\infty=0$ and under which $\alpha
_\infty>0$.
Theorem \ref{Thm-5} demonstrates that under mild conditions, $\alpha
_\infty$ of the $\FDR_L$ procedure
is lower than the counterpart of the $\FDR$ procedure.
These theoretical results demonstrate that the $\FDR_L$ procedure
alleviates the extent of the $\LIP$, thus
substantially facilitates the selection of user-specified control levels.
As observed from the middle and right panels of Figure \ref{Figure-2},
for control levels close to
zero,
the $\FDR_L$ procedure combined with either Method I or Method II
identifies a larger number of true significant locations than the
$\FDR$ procedure.
\end{itemize}

The rest of the paper is arranged as follows.
Section \ref{sec-2} reviews the conventional $\FDR$ procedure and introduces
$\alpha_\infty$
to characterize the $\LIP$.
Section \ref{sec-3} describes the proposed $\FDR_L$ procedure.
Its theoretical properties are established in Section \ref{sec-4}, where
Section \ref{sec-4.2}
explores the extent to which the $\FDR_L$ procedure alleviates $\LIP$.
Sections \ref{sec-5} and \ref{sec-6} present simulation comparisons of
the $\FDR$ and $\FDR_L$ procedures in 2D and 3D
dependent data, respectively.
Section \ref{sec-7} illustrates the computational simplicity and
detection effectiveness of the proposed
method for a real brain fMRI dataset for detecting the regions of activation.
Section \ref{sec-8} ends the paper with a brief discussion.
Technical conditions and detailed proofs are deferred to the \hyperref[app]{Appendix}.

\section{$\FDR$ and lack of identification phenomenon} \label{sec-2}

\subsection{Conventional $\FDR$ procedure}
We begin with a brief overview of the conventional $\FDR$ procedure
that is of particular relevance to the
discussion in Sections \ref{sec-3} and \ref{sec-4}.
For testing a family of null hypotheses,
$\{H_{0}(i)\}_{i=1}^n$,
suppose that $p_i$ is the $p$-value of the $i$th test.
Table \ref{Table-1} summarizes the outcomes.

%
%
\begin{table}[b]
\tablewidth=280pt
\caption{Outcomes from testing $n$ (null) hypotheses
$H_{0}(i)$ based on a significance rule}
\label{Table-1}
\begin{tabular*}{\tablewidth}{@{\extracolsep{\fill}}lccc@{}}
\hline
& $\bolds{H_{0}(i)}$ \textbf{retained} & $\bolds{H_{0}(i)}$ \textbf{rejected} & \textbf{Total} \\
\hline
$H_{0}(i)$ true & $U$ & $V$ & $n_0$ \\
$H_{0}(i)$ false & $T$ & $S$ & $n_1$ \\
[3pt]
Total & $W$ & $R$ & $n$ \\
\hline
\end{tabular*}
\end{table}

\citet{BH95} proposed a procedure that guarantees the
False Discovery Rate ($\FDR$)
to be less than or equal to a pre-selected value. Here, the $\FDR$ is
the expected
ratio of
the number of incorrectly rejected hypotheses to the total number of
rejected hypotheses with
the ratio defined to be zero if no hypothesis is rejected, that is,
$\FDR= E(\frac{V}{R\vee1})$
where $R\vee1 = \max(R,1)$.
A comprehensive overview of the development of the research in the area
of multiple testing can be found in \citet{BY01},
\citet{GW02}, \citet{S02}, \citet{DSB03},
\citet{E04},
\citet{STS04}, \citet{GW04},
\citet{LR05}, \citet{LRS05},
\citet{GRW06},
\citet{S06}, \citet{BH07} and \citet{W08}, among others.
Fan, Hall and Yao (\citeyear{FHY07}) addressed the issue on the number of
hypotheses that can be simultaneously tested when the $p$-values are
computed based on asymptotic approximations.

\citet{STS04} gave an empirical process definition of $\FDR$, by
%
%
\begin{equation} \label{b.1}
\FDR(t)=E\biggl\{\frac{V(t)}{R(t)\vee1}\biggr\},
\end{equation}
where $t$ stands for a threshold for $p$-values.
For realistic applications, \citet{S02} proposed
the point estimate of $\FDR(t)$ by
%
%
\begin{equation}\label{b.2}
\widehat\FDR(t) = \frac{W(\lambda)t}{\{R(t)\vee1\}(1-\lambda)},
\end{equation}
where $\lambda\in(0,1)$ is a tuning constant, and $W(t)$ is the
number of
nonrejections with a threshold $t$. The intuition of this will be
explained in Section \ref{sec-3.4}.
The pointwise limit of $\widehat\FDR(t)$ under assumptions (7)--(9) of
\citet{STS04}
is
%
%
\begin{equation}\label{b.3}
\widehat{\FDR}^{\infty}(t)=
\frac{[\pi_0\{1-G_0(\lambda)\}+\pi_1\{1-G_1(\lambda)\}] t}{\{\pi
_0 G_0(t)+\pi_1 G_1(t)\} (1-\lambda)},
\end{equation}
where $\pi_0 = \lim_{n\to\infty} n_0 / n$, $\pi_1=1-\pi_0$, and
$\lim_{n\to\infty} V(t)/n_0=G_0(t)$
and
$\lim_{n\to\infty} S(t)/n_1=G_1(t)$ are assumed to exist almost
surely for each $t\in(0, 1]$.
For a pre-chosen
level $\alpha$, a data-driven threshold for $p$-values is determined by
%
%
\begin{equation} \label{b.4}
t_{\alpha}(\widehat\FDR) = \sup\{0\le t \le1\dvtx \widehat\FDR(t) \le
\alpha\}.
\end{equation}
A null hypothesis is rejected if the corresponding $p$-value
is less than or equal to the threshold $t_{\alpha}(\widehat\FDR)$.
Methods (\ref{b.2}) and (\ref{b.4}) form the basis for the
conventional $\FDR$ procedure.

\subsection{Proposed measure for lack of identification phenomenon}
Recall that the $\FDR$ procedure is essentially a threshold-based
approach for multiple testing
problems, where the data-driven threshold $t_{\alpha}(\widehat\FDR)$
plays a key role. It is clearly seen from
(\ref{b.4})
that $t_{\alpha}(\widehat\FDR)$ hinges on both the estimates $\widehat
\FDR(t)$
devised, as well as the control level $\alpha$ specified.

Using
(\ref{b.2}),
we observe that the corresponding $t_{\alpha}(\widehat\FDR)$ is a
nondecreasing function of
$\alpha$.
This indicates that
for the $\FDR$ procedure, as $\alpha$ decreases below
$\inf_{0 < t\le1}\widehat\FDR(t)$,
the threshold $t_{\alpha}(\widehat\FDR)$ will drop to zero and
accordingly,
the $\FDR$ procedure {can only reject those hypotheses with $p$-values
exactly equal to zero}.
We call this phenomenon ``\textit{lack of identification}.''

To better quantify the ``\textit{lack of identification phenomenon}''
($\LIP$), the limiting forms
of $\widehat\FDR(t)$ as $n\to\infty$ will be examined.
\begin{definition} \label{def-1}
For estimation methods $\widehat\FDR(t)$ in (\ref{b.2}),
define
\[
\alpha_{\infty}^{\FDR} = \inf_{0< t\le1} \widehat\FDR^\infty(t),
\]
where
$\widehat\FDR^\infty(t)$
is defined in (\ref{b.3}).
Define the
endurance by
$E_{\FDR} = 1 - \alpha_{\infty}^{\FDR}$.
\end{definition}

Notice that the existence of $\alpha_{\infty}^{\FDR}>0$ implies the
occurrence of the $\LIP$:
in real data applications with a moderately large number $n$ of hypotheses,
the $\FDR$ procedure loses the identification capability when the
control level $\alpha$ is
close to or smaller than $\alpha_{\infty}^{\FDR}$.
On the other hand,
the case $\alpha_{\infty}^{\FDR}=0$ rules out the possibility of the
$\LIP$.
Henceforth, the smaller the $\alpha_{\infty}^{\FDR}$, the higher
endurance of the
corresponding $\widehat\FDR$, and the less likely the $\LIP$ happens.
In other words, an $\FDR$ estimation approach
with a higher endurance is more capable of adopting a smaller control
level, thus
reducing the extent of the $\LIP$ problem.
We will revisit this issue in Section \ref{sec-4.2} after introducing
the proposed $\FDR_L$ procedure.

\section{Proposed $\FDR_L$ procedure for imaging data} \label{sec-3}

Consider a set of spatial signals $\{\mu(v)\dvtx v \in\calV\subseteq\ZZ
^{d}\}$
in a 2D plane ($d=2$) or a 3D space ($d=3$), where
$\mu(v) = 0$ for $v\in\calV_0$,
$\mu(v) \ne0$ for $v\in\calV_1$
and $\calV_0 \cup\calV_1=\calV$.
Here $\calV_0$ and $\calV_1$ are unknown sets.
A common approach for detecting the presence of the spatial signals
consists of two stages. In the first stage,
test the hypothesis
\[
H_0(v)\dvtx \mu(v) = 0 \quad\versus\quad H_1(v)\dvtx \mu(v) \ne0
\]
at each location $v$. The corresponding $p$-value is denoted by
$p(v)$. In the second stage, a multiple testing procedure, such as the
conventional
$\FDR$ procedure, is applied to
the collection, $\{p(v)\dvtx v \in\calV\subseteq\ZZ^{d}\}$, of $p$-values.

In the second stage, instead of using the original $p$-value, $p(v)$,
at each $v$, we propose to use a
local aggregation of $p$-values at points located adjacent to $v$. We
summarize the procedure as follows.
\begin{step}\label{step1}
Choose a local neighborhood with size $k$.
\end{step}
\begin{step}\label{step2}
At each grid point $v$, find the set $N_v$ of
its neighborhood points,
and the set $\{p(v')\dvtx v'\in N_v\}$ of the corresponding $p$-values.
\end{step}
\begin{step}\label{step3}
At each grid point $v$, apply a
transformation $f\dvtx [0, 1]^{k} \mapsto[0, 1]$ to
the set of $p$-values
in Step \ref{step2}, leading to a ``locally aggregated'' quantity,
$p^*(v)=f(\{p(v')\dvtx v'\in N_v\})$.
\end{step}
\begin{step}\label{step4}
Determine a data-driven threshold
for $\{p^*(v)\dvtx v\in\calV\subseteq\ZZ^{d}\}$.
\end{step}

For
notational clarity, we denote by $\{p_i^*\}_{i=1}^n$ the collection of
``locally aggregated'' $p^*$-values,
$\{p^*(v)\dvtx v\in\calV\subseteq\ZZ^{d}\}$.
Likewise, the notation $U^*(t)$, $V^*(t)$, $T^*(t)$, $S^*(t)$, $W^*(t)$
and $R^*(t)$ can be defined as
in Section \ref{sec-2}, with $p_i$ replaced by $p_i^*$.
For instance, $V^*(t) = \sum_{i=1}^n \ID\{H_0(i)$ is true, and $p_i^*
\le t\}$
and
$R^*(t) = \sum_{i=1}^n \ID(p_i^* \le t)$,
with $\ID(\cdot)$ an indicator function.
Accordingly, the false discovery rate based on utilizing the locally
aggregated $p_i^*$-values becomes
%
%
\begin{equation} \label{c.1}
\FDR_L(t)=E\biggl\{\frac{V^*(t)}{R^*(t)\vee1}\biggr\}.
\end{equation}
As a comparison, $\FDR(t)$ in (\ref{b.1}) corresponds to the use of
the original $p$-values.

\subsection{Choice of neighborhood and choice of $f$}
As in \citet{RS00}, the set of neighbors for each data point
can be assigned in a variety of ways,
by choosing the $k$ nearest neighbors in Euclidean distance,
by considering all data points within a ball of fixed radius
or by using some prior knowledge.

For the choice of the transformation function, $f$, one candidate is
the median filter, applied to the neighborhood $p$-values, without
having to specify
particular forms of spatial structure.
A discussion on other options for $f$ can be found in Section \ref{sec-8}.
Unless otherwise stated, this paper focuses on the median filtering.

\subsection{Statistical inference for $p^*$-values: Method $\mathrm{I}$} \label{sec-3.2}

Let $G^*(\cdot)$ be the cumulative distribution function of a
``locally aggregated'' $p^*$-value corresponding to the true null hypothesis.
Let $\widetilde G^*(\cdot)$ be the sample distribution of $\{p^*(v):
v\in
\calV_0\}$.
Recall that the original $p$-value corresponding to the true null
hypothesis is uniformly distributed on
the interval $(0,1)$. In contrast,
the distribution $G^*(\cdot)$ for a ``locally aggregated'' $p^*$-value
is typically nonuniform.
This indicates that
a significance rule based on $p$-values is not directly applicable to
the significance rule based on $p^*$-values.
For the median\vspace*{1pt} operation $f$, we propose two methods for estimating
$\widetilde G^*(\cdot)$.
Method I is particularly useful for large-scale imaging datasets,
whereas Method II
is useful for data of limited resolution.

Method I is motivated from the observation: if the original
$p$-values are independent and uniformly
distributed on the interval $(0,1)$, then the median aggregated
$p^*$-value follows a
Beta distribution. More precisely,
if the neighborhood size $k$ is an odd integer, then the median
aggregated $p^*$-value conforms to the
%
%
\begin{equation} \label{c.2}
\Betadist\bigl((k+1)/2, (k+1)/2\bigr)
\end{equation}
distribution [\citet{CB90}].
If $k$ is an even integer,
the median aggregated $p^*$-value is distributed as a random variable
$(X+Y)/{2}$, where
$(X,Y)$ has the joint probability density function
${k!}/{\{(k/2-1)!\}^2} x^{k/2-1} (1-y)^{k/2-1} \ID(0<x<y<1)$.
Thus, as long as the resolution of the experiment data and imaging
technique keeps
improving, so that the proportion of boundary grid points
(corresponding to those with neighborhood
intersected with both $\calV_0$ and $\calV_1$) decreases and
eventually shrinks to zero,
$G^*(\cdot)$ will tend to the Beta distribution in (\ref{c.2}).

Following this argument, if the original $p$-values corresponding to
the true null hypotheses are independent
and uniformly distributed [see, e.g., \citet{Vaart98}, page 305],
the median aggregated $p^*$-values corresponding to the true null
hypotheses will approximately be
symmetrically distributed about $0.5$.
Thus, assuming that the number of false null hypothesis with $p_i^* >
0.5$ is negligible, the total number of true null hypotheses, $n_0$,
is approximately
$2\sum_{i=1}^n\ID(p_i^*>0.5)+\sum_{i=1}^n\ID(p_i^*=0.5)$, and the number
of true null hypotheses with $p^*$-values smaller than or equal to $t$
could be estimated by
$\sum_{i=1}^n\ID\{p_i^* \ge(1-t)\}$, for small values of $t$. Here,
owing to the symmetry, we use the upper tail to compute the proportion
to mitigate the bias caused by the data from the alternative
hypotheses. Hence,
$\widetilde G^*(t)$ can be estimated by the empirical distribution function,
%
%
\begin{equation} \label{c.3}
\widehat G^*(t) =
\cases{
\dfrac{\sum_{i=1}^n \ID\{p_i^* \ge(1-t)\}}{2\sum_{i=1}^n\ID
(p_i^*>0.5)+\sum_{i=1}^n\ID(p_i^*=0.5)},\vspace*{3pt}\cr
\qquad \mbox{if $0 \le t\le0.5$}, \vspace*{3pt}\cr
1- \dfrac{\sum_{i=1}^n \ID(p_i^*>t)}{2\sum_{i=1}^n\ID
(p_i^*>0.5)+\sum_{i=1}^n\ID(p_i^*=0.5)},\vspace*{3pt}\cr
\qquad \mbox{if $0.5 <t \le1$}.}
\end{equation}
A modification
of the Glivenko--Cantelli theorem shows that ${\sup_{0\le t\le1}}
|\widehat
G^*(t)-G^*(t)| = o(1)$ almost surely as $n\to\infty$.
This method is distribution free,
computationally fast and applicable when the $p^*$-values under the
null hypotheses are not too skewedly distributed.

An alternative approach for approximating $\widetilde G^*(\cdot)$ is
inspired by the central limit
theorem. If the neighborhood size $k$ is
reasonably large (e.g., $k\ge5$ if the original $p$-values
corresponding to the true null hypotheses are independent and uniformly
distributed), then $\widetilde G^*(\cdot)$
could be approximated by a normal distribution centered at $0.5$.
This normal approximation scheme may be exploited in the situation
(which rarely occurs, though)
when the original $p$-values corresponding to the true null hypotheses
are independent but asymmetric about $0.5$ (when the null distribution
function of the test statistic
is discontinuous).

\subsection{Refined method for estimating $\widetilde G^*(\cdot)$: Method
$\mathrm{II}$} \label{sec-3.3}

More generally, we consider spatial image data of limited resolution.
Recall the neighborhood size $k$ of a voxel $v$ in the paper includes
one for
$v$ itself.
Let $\mathsf{n}_1(v)$ denote the number of points in $N_v$ that belong
to $\calV_1$.
Thus for any grid point $v\in\calV_0$, $\mathsf{n}_1(v)$ takes
values $\{0,1,\ldots,k-1\}$.
Set
\[
\theta_{n,j}=P\{\mathsf{n}_1(v) = j\},\qquad
Q_j^*(t)=P\{p^*(v) \le t |\mathsf{n}_1(v) = j\}.
\]
Clearly, $\sum_{j=0}^{k-1}\theta_{n,j} = 1$. Therefore, the C.D.F. of
$p^*(v)$ for a grid point $v\in\calV_0$ is given by
%
%
\begin{equation} \label{c.4}
G^*(t)
=\theta_{n,0} Q_0^*(t) + \theta_{n,1} Q_1^*(t) + \cdots+ \theta
_{n,k-1}Q_{k-1}^*(t),
\end{equation}
where $Q_0^*(t)$ corresponds to, for independent tests, the Beta
distribution function in (\ref{c.2}).

Likewise, we obtain
\[
\widetilde G^*(t) = \sum_{j=0}^{k-1}\widetilde\theta_{n,j} \widetilde Q_j^*(t),
\]
where
$\widetilde\theta_{n,j} = {\# \calV_0^{(j)}}/{n_0}$ is the proportion of
$v \in\calV_0$ with $j$ neighboring grid points in $\calV_1$,
and $\widetilde Q_j^*(t) = {\sum_{v \in\calV_0^{(j)}} \ID\{p^*(v) \le
t\}}/{\# \calV_0^{(j)}}$ is the sample\vspace*{-2pt} distribution of $\{p^*(v): v
\in\calV_0^{(j)}\}$,
with $\# A$ denoting the number of elements in a set $A$ and
$\calV_0^{(j)} = \{v \in\calV_0\dvtx \mathsf{n}_1(v)=j\}$. Clearly, if
the original $p$-values corresponding to the true null hypotheses are
block dependent, then, by the Glivenko--Cantelli theorem, $\sup_{0\le
t \le1} |\widetilde G^*(t) - G^*(t)| = o(1)$ almost surely, as $n\to
\infty$.

We propose the following
Method II
to estimate $\widetilde G^*(t)$:
\begin{enumerate}
\item Obtain estimates $\widehat n_0$ and $\widehat n_1 = n-\widehat
n_0$ of $n_0$
and $n_1$, respectively. One possible estimator of $n_0$ is
$\widehat n_0 = {\sum_{i=1}^n \ID(p_i^* > \lambda)}/\{1-\widehat
G^*(\lambda
)\}$,
for some tuning parameter $\lambda$.

\item Define
$\widehat\calV_1 = \{v\in\calV\dvtx p^*(v) \le p_{(\widehat n_1)}^*\}$, where
$\{p_{(i)}^*\}_{i=1}^n$ denote the order statistics of
$\{p_i^*\}_{i=1}^n$.
Define $\widehat\calV_0 = \{v\in\calV\dvtx p^*(v) > p_{(\widehat n_1)}^*\}$.

\item Set $\widehat\calV_0^{(j)}=\{v \in\widehat\calV_0\dvtx \mathsf
{n}_1(v)=j\}$. Estimate $\widetilde\theta_{n,j}$, $j=0,\ldots,k-1$, by
$\widehat\theta_{n,j} = {\# \widehat\calV_0^{(j)}}/{\widehat n_0}$.

\item
For $j = 0$, estimate $\widetilde Q_0^*(t)$ by
$\widehat Q_0^*(t) = \widehat G^*(t)$,
the estimator of $\widetilde G^*(t)$ by Method~I in Section \ref{sec-3.2}.
To estimate $\widetilde Q_j^*(t)$, $j=1,\ldots, k-1$, for each $v \in
\widehat
\calV_0^{(0)}$, collect its neighborhood $p$-values, randomly exclude $j$
of them and obtain the set
$D_j(v)$ for the remaining neighborhood $p$-values. Randomly sample $j$
grid points from $\widehat\calV_1$ and
collect their corresponding $p$-values in a set $A_j(v)$. Compute
the median, $\widehat p_j^*(v)$, of $p$-values in $D_j(v)\cup A_j(v)$.
Estimate $\widetilde Q_j^*(t)$ by
$\widehat Q_j^*(t) = {\sum_{v\in\widehat\calV_0^{(0)}}
\ID\{\widehat p_j^*(v)\le t\}}/{\#\widehat\calV_0^{(0)}}$.\vadjust{\goodbreak}

\item Combining (\ref{c.4}),
$\widetilde G^*(t)$ is estimated by
$\widehat G_c^*(t) =
\sum_{j=0}^{k-1}\widehat\theta_{n,j} \widehat Q_j^*(t)$.
\end{enumerate}

\subsection{Significance rule for $p^*$-values} \label{sec-3.4}
Using the locally aggregated $p^*$-values, we can estimate
$\FDR_L(t)$ defined in (\ref{c.1}) by either
%
%
\begin{equation} \label{c.9}
\widehat\FDR_L(t) = \frac{W^*(\lambda)\widehat G^*(t)}{\{R^*(t)\vee1\}\{
1-\widehat
G^*(\lambda)\}},
\end{equation}
using Method I, or
%
%
\begin{equation} \label{c.10}
\widehat\FDR_L(t) = \frac{W^*(\lambda)\widehat G_c^*(t)}{\{R^*(t)\vee1\}
\{
1-\widehat
G_c^*(\lambda)\}},
\end{equation}
using Method II. The logic behind this estimate is the following.
If we
choose $\lambda$ far enough from zero, then the number of
nonrejections, $W^*(\lambda)$, is roughly $U^*(\lambda)$.
Using this, we have
\[
V^*(\lambda) \approx n_0 \widetilde G^*(\lambda) \approx
\{V^*(\lambda) + W^*(\lambda)\} \widetilde G^*(\lambda).
\]
Solving the above equation suggests an estimate of $V^*(\lambda)$ by
${W^*(\lambda)\widetilde G^*(\lambda)}/\{1-\widetilde G^*(\lambda)\}$.
Now, using
$V^*(t) / V^*(\lambda) \approx\widetilde{G}^*(t)/\widetilde{G}^*(\lambda
)$, we obtain that at a threshold~$t$, $V^*(t)$ can be estimated by
${W^*(\lambda)\widetilde G^*(t)}/\{1-\widetilde G^*(\lambda)\}$.
This together with the definition of $\FDR_L(t)$ in (\ref{c.1})
suggests the estimate\vspace*{1pt} in
(\ref{c.9}). Interestingly,
in the particular case of $p_i^*\equiv p_i$ and $\widehat G^*(t)=t$ [or
$\widehat G_c^*(t)=t$],
$\widehat\FDR_L(t)$ coincides with $\widehat\FDR(t)$ defined in (\ref{b.2}).

For a given control level $\alpha$, a null hypothesis is rejected if
the associated
$p^*$-value is smaller than or equal to the threshold,
%
%
\begin{equation} \label{c.11}
t_{\alpha}(\widehat\FDR_L) \equiv\sup\{0\le t \le1\dvtx
\widehat\FDR_L(t)\le\alpha\}.
\end{equation}
This data-driven threshold for $p^*$-values together with the point
estimation method (\ref{c.9})
[or (\ref{c.10})] for the false discovery rates comprises the proposed
$\FDR_L$ procedure.

\section{Properties of the $\FDR_L$ procedure} \label{sec-4}

\subsection{Asymptotic behavior} \label{sec-4.1}
This section explores the asymptotic behavior of the $\FDR_L$
procedure under weak dependence of $p$-values.
Technical assumptions are given in Condition \ref{CondA} in the \hyperref[app]{Appendix},
where Conditions \hyperlink{CondA1}{A1}--\hyperlink{CondA3}{A3} are similar to
assumptions $(7)$--$(9)$ of \citet{STS04}. Thus the type of dependence
in Condition \hyperlink{CondA2}{A2} includes
finite block dependence,
and certain mixing dependence.
Theorems \ref{Thm-1}--\ref{Thm-3} can be considered a generalization
of \citet{STS04}
from a single $p$-value to locally aggregating a number $k$ of
$p$-values with $k>1$.\vadjust{\goodbreak}

Theorem \ref{Thm-1} below reveals that the proposed estimator $\widehat
\FDR_L(t)$
controls the $\FDR_L(t)$ simultaneously for all $t\ge\delta$ with
$\delta>0$, and in turn supplies a
conservative estimate of $\FDR_L(t)$.
\begin{theorem} \label{Thm-1}
Assume Condition \ref{CondA} in Appendix \ref{appA}. For each $\delta> 0$,
\[
\lim_{n\to\infty} \inf_{t\ge\delta} \biggl\{\widehat\FDR_L(t)-\frac
{V^*(t)}{R^*(t)\vee1}\biggr\}
\ge0
\]
and
\[
\lim_{n\to\infty} \inf_{t\ge\delta} \{\widehat\FDR_L(t)-\FDR
_L(t)\}
\ge0
\]
with probability one.
\end{theorem}

To show that the proposed $\widehat\FDR_L(t)$ asymptotically provides a strong
control of $\FDR_L(t)$, we define
%
%
\begin{equation} \label{d.1}
\widehat\FDR_L^{\infty}(t) =
\frac{[\pi_0\{1-G_0^*(\lambda)\}+\pi_1\{1-G_1^*(\lambda)\}]
G^{*\infty}(t)}{\{\pi_0 G_0^*(t)+\pi_1 G_1^*(t)\} \{1- G^{*\infty
}(\lambda)\}},
\end{equation}
which is the pointwise limit of $\widehat\FDR_L(t)$ under Condition
\ref{CondA} in Appendix \ref{appA},
where it is assumed that $\pi_0 = \lim_{n \to\infty} n_0/n$, and
$\lim_{n \to\infty} {V^*(t)}/{n_0} = G_0^*(t)$ and
$\lim_{n \to\infty} {S^*(t)}/{n_1} = G_1^*(t)$ exist almost surely
for each $t\in(0,1]$,
and $G^{*\infty}(t) = \lim_{n\to\infty} G^*(t)$.
\begin{theorem} \label{Thm-2}
Assume Condition \ref{CondA} in Appendix \ref{appA}. If there is a $t\in(0,1]$
such that $\widehat\FDR_L^{\infty}(t)<\alpha$, then
$
\limsup_{n\to\infty} \FDR_L(t_{\alpha}(\widehat\FDR_L)) \le
\alpha$.
\end{theorem}

Theorem \ref{Thm-3} states that the random thresholding rule
$t_{\alpha}(\widehat\FDR_L)$ converges
to the deterministic rule $t_{\alpha}(\widehat\FDR_L^{\infty})$.
\begin{theorem} \label{Thm-3}
Assume Condition \ref{CondA} in\vspace*{1pt} Appendix \ref{appA}. If
$\widehat\FDR_L^{\infty}(\cdot)$ has a nonzero derivative at the point
$t_{\alpha}(\widehat\FDR_L^{\infty}) \in(0, 1)$, then
$\lim_{n\to\infty} t_{\alpha}(\widehat\FDR_L) = t_{\alpha}(\widehat\FDR
_L^{\infty})$
holds almost surely.
\end{theorem}

\subsection{Conditions for lack of identification phenomenon} \label{sec-4.2}
\begin{definition} \label{def-2}
For estimation methods $\widehat\FDR_L(t)$ in (\ref{c.9}) [or (\ref{c.10})],
define
\[
\alpha_{\infty}^{\FDR_L} = \inf_{0< t\le1} \widehat\FDR_L^\infty(t),
\]
where
$\widehat\FDR_L^\infty(t)$
is defined in
(\ref{d.1}).
\end{definition}

Theorem \ref{Thm-4} establishes conditions under which the $\LIP$
does or does not take place with the $\FDR$ and $\FDR_L$ procedures.
It will be seen that the conditions are characterized by the null and
alternative
distributions of the test statistics, without relying on
the configuration of the neighborhood used in the $\FDR_L$ procedure.
Theorem \ref{Thm-5} demonstrates that $\alpha_{\infty}^{\FDR} \ge
\alpha_{\infty}^{\FDR_L}$ under mild conditions,
thus the $\FDR_L$ procedure reduces the extent of the $\LIP$.
For expository brevity, we assume the test statistics
are independent,
which can be relaxed.
\begin{theorem} \label{Thm-4}
Let $\{T(v)\dvtx v\in\calV\subseteq\ZZ^{d}\}$ be the set of test statistics
for testing the presence of the spatial signals $\{\mu(v)\dvtx v\in\calV
\subseteq\ZZ^{d}\}$. Consider the one-sided testing problem,
%
%
\begin{equation} \label{d.2}
H_0(v)\dvtx \mu(v)=0 \quad\versus\quad H_1(v)\dvtx \mu(v)>0.
\end{equation}
For $j=0$ and $j=1$, respectively,
assume that $T(v)$, corresponding to the true $H_j(v)$, are i.i.d.
random variables having
a cumulative distribution function $F_j$
with
a probability density function $f_j$.
Assume that the neighborhood size $k\ge3$ used in the $\FDR_L$
procedure is an odd integer
and that the proportion of boundary grid points within $\calV_0$
shrinks to zero, as $n\to\infty$,
that is, \mbox{$\lim_{n\to\infty} {\# \calV_0^{(0)}}/{n_0} =1$}, where
$
\calV_0^{(0)} = \{v \in\calV\dvtx \mu(v') = 0$ for any $v' \in N_v \}.
$
Assume Condition \textup{\hyperlink{CondA1}{A1}} in Appendix \ref{appA}.
Let $x_0=F_0^{-1}(1) = \inf\{t\dvtx F_0(t) = 1\}
$.
\begin{longlist}[$\mathrm{II}.$]
\item[$\mathrm{I}.$] If
$
\lim_{x\to x_0-} \frac{f_1(x)}{f_0(x)} = \infty,
$
then $\alpha_{\infty}^{\FDR} = 0$ and $\alpha_{\infty}^{\FDR_L} = 0$.

\item[$\mathrm{II}.$] If
$
\limsup_{x\to x_0-}\frac{f_1(x)}{f_0(x)} < \infty,
$
then $\alpha_{\infty}^{\FDR} >0$ and $\alpha_{\infty}^{\FDR_L} >0$.
\end{longlist}
\end{theorem}
\begin{theorem} \label{Thm-5}
Assume the conditions\vspace*{1pt} in Theorem \ref{Thm-4}.
Suppose that $f_0(\cdot)$ is supported in an interval;
$f_1(x)\le f_0(x)$ for any $x\le F_0^{-1}(0.5)$;
$1-F_0(F_1^{-1}(0.5))\le\lambda\le0.5$.
Then
$\alpha_{\infty}^{\FDR} \ge\alpha_{\infty}^{\FDR_L}$.
\end{theorem}

Corollaries \ref{corollary-1} and \ref{corollary-2} below provide
concrete applications of Theorems \ref{Thm-4} and~\ref{Thm-5}.
The detailed verifications are omitted.
\begin{corollary} \label{corollary-1}
Assume the conditions in Theorem \ref{Thm-4}.
Suppose that
the distribution $F_0$ is $N(0, 1)$ and
the distribution $F_1$ is $N(C, \sigma^2)$,
where $\sigma\in(0, \infty)$ and $C \in(0, \infty)$ are constants.
\begin{longlist}[$\mathrm{II}.$]
\item[$\mathrm{I}.$] If $\sigma\ge1$, then $\alpha_{\infty}^{\FDR
} = 0$ and $\alpha_{\infty}^{\FDR_L} = 0$.

\item[$\mathrm{II}.$] If $0 < \sigma< 1$, then $\alpha_{\infty
}^{\FDR} > 0$ and $\alpha_{\infty}^{\FDR_L} > 0$.
Moreover,
if $\exp\{-(C/\sigma)^2/2\} /\break{\sigma} \le1$ and
$1- F_0(C)\le\lambda\le0.5$,
then $\alpha_{\infty}^{\FDR} \ge\alpha_{\infty}^{\FDR_L}$.
\end{longlist}
\end{corollary}
\begin{corollary} \label{corollary-2}
Assume the conditions in Theorem \ref{Thm-4}.
Suppose that
the distribution $F_0$ is that of a Student's $t_{d_0}$ variate with
$d_0$ degrees of freedom and
the distribution $F_1$ is that of $C$ plus a Student's $t_{d_1}$
variate with $d_1$ degrees of freedom,
where $C \in(0, \infty)$ is a constant.
\begin{longlist}[$\mathrm{II}.$]
\item[$\mathrm{I}.$] If $d_0 > d_1$, then $\alpha_{\infty}^{\FDR}
= 0$ and $\alpha_{\infty}^{\FDR_L} = 0$.\vadjust{\goodbreak}

\item[$\mathrm{II}.$] If $1 \le d_0 \le d_1$, then $\alpha_{\infty
}^{\FDR} > 0$ and $\alpha_{\infty}^{\FDR_L} > 0$.
Moreover,
if $d_0=d_1$ and
$1-F_0(C)\le\lambda\le0.5$,
then $\alpha_{\infty}^{\FDR} \ge\alpha_{\infty}^{\FDR_L}$.
\end{longlist}
\end{corollary}
\begin{remark}
For illustrative simplicity, a one-sided testing problem (\ref{d.2})
is focused upon.
Two-sided testing problems can similarly be treated and we omit the details.
\end{remark}

\subsection{An illustrative example of $\alpha_{\infty}^{\FDR} >
\alpha_{\infty}^{\FDR_L}>0$} \label{sec-4.3}

\begin{figure}[b]

\includegraphics{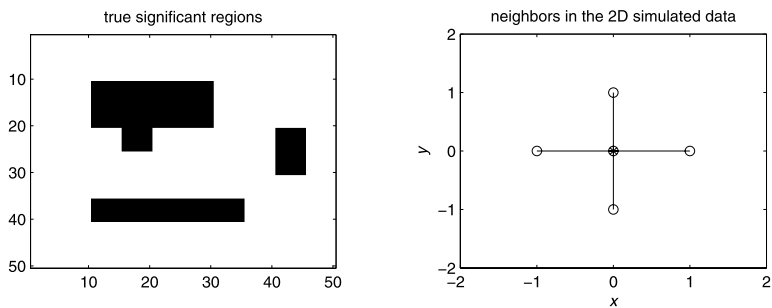}

\caption{Left panel: the true significant regions for the $\mathrm
{2D}$ simulated data sets.
Right panel: neighbors of a point at $(x,y)$ used in the $\FDR_L$
procedure for
$\mathrm{2D}$ simulated data.}
\label{Figure-1}
\end{figure}

Consider a pixelated 2D image dataset
consisting of $n = 50\times50$
pixels, illustrated in the left panel of Figure \ref{Figure-1},
where
the black rectangles represent the true significant regions $\calV_1$
with $n_1 = 0.16 \times n$ pixels
and
the white background serves as the true nonsignificant regions $\calV
_0$ with $n_0 = n-n_1$ pixels.
The data are simulated from the model,
\[
Y(i,j)=\mu(i,j)+\varepsilon(i,j),\qquad i,j=1,\ldots,50,
\]
where
the signals are
$\mu(i,j)=0$ for $(i,j)\in\calV_0$, and
$\mu(i,j)=C$ for $(i,j)\in\calV_1$ with a constant $C\in(0,\infty)$,
and
the error terms $\{\varepsilon(i,j)\}$ are i.i.d. following the centered
$\operatorname{Exp}(1)$ distribution.
At each site $(i,j)$, the observed data $Y(i,j)$ is the (shifted)
survival time and used as the test statistic for
testing $\mu(i,j)=0$ versus $\mu(i,j)>0$.
Clearly,
all test statistics corresponding to the true null hypotheses are i.i.d.
having the probability density function $f_0(x)=\exp\{-(x+1)\}\ID
(x+1>0)$; likewise,
all test statistics in accordance with the true alternative hypotheses
are i.i.d.
having the density function $f_1(x) = \exp\{-(x+1-C)\} \ID(x+1>C)$.
It is easily seen
that $x_0
= \infty$, and
$\limsup_{x\to\infty}{f_1(x)}/{f_0(x)} = \exp(C)<\infty$.\vspace*{1pt}
An appeal to Theorem \ref{Thm-4} yields $\alpha_{\infty}^{\FDR}>0$
and $\alpha_{\infty}^{\FDR_L}>0$, and thus
both the $\FDR$ and $\FDR_L$ procedures will encounter the $\LIP$.
Moreover, if
$C > \log(2)$, $\exp(-C)/2 \le\lambda\le0.5$
and the neighborhood size $k\ge3$ is an odd integer,\vspace*{1pt}
then sufficient conditions in
Theorem \ref{Thm-5} are satisfied and hence $\alpha_{\infty}^{\FDR}
\ge\alpha_{\infty}^{\FDR_L}$.

%
%

Actual computations indicate that in this example,
as long as $C>\log(4)$,
$\alpha_{\infty}^{\FDR_L}$ is considerably smaller than $\alpha
_{\infty}^{\FDR}$,
indicating that the
$\FDR_L$ procedure can adopt a control level
much smaller than that of the conventional $\FDR$ procedure without excessively
encountering the $\LIP$.
For example, set $\lambda= 0.1$; assume that
the neighborhood in the $\FDR_L$ procedure is depicted in the right
panel of Figure
\ref{Figure-1}, that is, $k=5$.
Table \ref{Table-2} compares values of $\alpha_{\infty}^{\FDR}$ and
$\alpha_{\infty}^{\FDR_L}$
for $C=\log(4j)$, $j=2,\ldots,9$.
Refer to (\ref{C.2}) and (\ref{C.5}) in Appendix \ref{appC} for detailed
derivations of $\alpha_\infty^{\FDR}$ and $\alpha_{\infty}^{\FDR
_L}$, respectively.

%
%
\begin{table}
\tabcolsep=4pt
\caption{Comparing $\protect\vphantom{{{a^{\sum}}^{\int}}^{\sum}}\alpha_{\infty}^{\FDR}$ and $\alpha
_{\infty}^{\FDR_L}$} \label{Table-2}
\begin{tabular*}{\tablewidth}{@{\extracolsep{\fill}}lcccccccc@{}}
\hline
$\bolds{C}$ & $\bolds{\log(8)}$ & $\bolds{\log(12)}$ & $\bolds{\log(16)}$ & $\bolds{\log(20)}$
& $\bolds{\log(24)}$ & $\bolds{\log(28)}$ & $\bolds{\log(32)}$ & $\bolds{\log(36)}$ \\
\hline
$\alpha_{\infty}^{\FDR}$ & 0.4130 & 0.3043 & 0.2471 & 0.2079 &
0.1795 & 0.1579 & 0.1409 & 0.1273 \\[4pt]
$\alpha_{\infty}^{\FDR_L}$ & 0.0103 & 0.0030 & 0.0013 & 0.0007 &
0.0004 & 0.0002 & 0.0002 & 0.0001 \\
\hline
\end{tabular*}
\end{table}

%
%
\begin{figure}[b]

\includegraphics{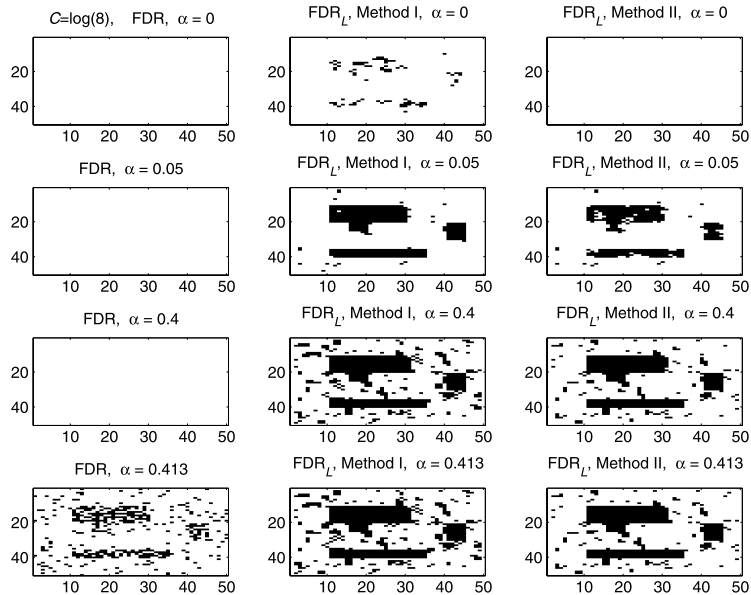}

\caption{Lack of identification phenomenon when $\alpha$
varies from $0$ to $\alpha_{\infty}^{\FDR} = 0.4130$. The sites that
are called statistically significant based on the realization are shown
in black. Left panels: the $\FDR$ procedure. Middle panels: the
$\FDR_L$ procedure using Method $\mathrm{I}$. Right panels: the
$\FDR_L$ procedure using Method $\mathrm{II}$.} \label{Figure-2}
\end{figure}

To better visualize the $\LIP$ from limited data,
Figure \ref{Figure-2} compares the regions detected
as significant by the $\FDR$
and $\FDR_L$ procedures for $C=\log(8)$ based on one realization of
the simulated data.
It is observed from Figure \ref{Figure-2} that
for $\alpha$ between $0$ and $0.4$, the $\FDR$ procedure lacks the
ability to detect statistical significance;
as $\alpha$ increases to $0.413$ (which is the limit $\alpha_{\infty
}^{\FDR}=0.413$ as calculated in Table \ref{Table-2})
and above, some significant results emerge. In contrast,
for $\alpha$ close to $0$,
both Method I and Method II for the $\FDR_L$ procedure are able to
deliver some significant results.
Similar plots to those in Figure \ref{Figure-2} are obtained with
other choices of $C$
and hence are omitted for lack of space.

\section{Simulation study: 2D dependent data} \label{sec-5}

\subsection{Example $1$} \label{sec-5.1}
\begin{figure}[t!]

\includegraphics{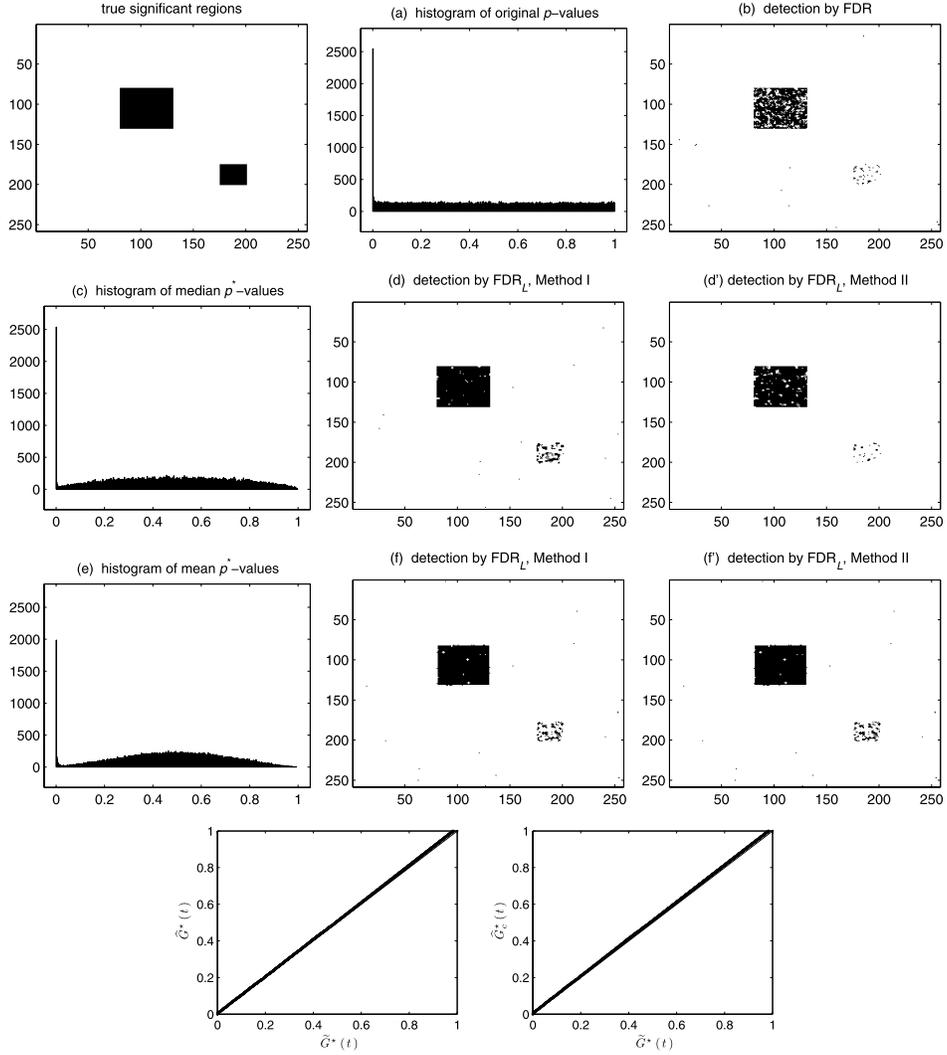}

\caption{Comparison of the $\FDR$ and $\FDR_L$ procedures
for Example $1$. In the first row, left: true significant regions shown
in black; middle: histogram of the original $p$-values; right:
significant regions detected by the $\FDR$ procedure. In the second
row, left: histogram of the $p^*$-values using the median filter;
middle and right: significant regions detected by the $\FDR_L$
procedure using Methods~$\mathrm{I}$ and $\mathrm{II}$, respectively.
In the third row, left: histogram of the $p^*$-values using the mean
filter; middle and right: significant regions detected by the $\FDR_L$
procedure\vspace*{1pt} using Methods $\mathrm{I}$ and $\mathrm{II}$, respectively.
In the bottom row, left: $\widehat G^*(t)$ versus $\widetilde G^*(t)$;
right: $\widehat G_c^*(t)$ versus $\widetilde G^*(t)$; straight line:
the $45$
degree reference line. Here $\alpha=0.01$ and $\lambda=0.1$.}
\label{Figure-3}
\vspace*{12pt}
\end{figure}

To illustrate the distinction between the $\FDR_L$ and the
conventional $\FDR$ procedures,
we present simulation studies.
The true significant regions are displayed as two black rectangles in
the top left panel of Figure \ref{Figure-3}.
The data are generated according to the model
%
%
\begin{equation} \label{e.1}
Y(i,j)=\mu(i,j)+\varepsilon(i,j),\qquad i,j=1,\ldots,258,
\end{equation}
where
the signals are
$\mu(i,j)=0$ for $(i,j)\in\calV_0$,
$\mu(i,j)=4$ in the larger black rectangle and $\mu(i,j)=2$ in the
smaller black rectangle.
The errors $\{\varepsilon(i,j)\}$ have zero-mean, unit-variance and are
\textit{spatially dependent},\vspace*{1pt} by taking
$\varepsilon(i,j)=
\{e(i-1,j) +e(i, j) +e(i+1,j) +e(i,
j-1)+e(i,j+1)\}/\sqrt{5}$,\vspace*{1pt}
where $\{e(i,j)\}_{i,j=0}^{259}$ are i.i.d. $N(0,1)$.
At each pixel $(i,j)$, $Y(i,j)$ is used as the test statistic for
testing $\mu(i,j)=0$ against $\mu(i,j) > 0$.

Both $\FDR$ and $\FDR_L$ procedures are preformed at a common control
level $0.01$, with the tuning
constant $\lambda= 0.1$.
In the $\FDR_L$ procedure, the neighborhood of a point at $(x,y)$ is
taken as in the right panel of Figure \ref{Figure-1}.
The histogram of the original $p$-values plotted in Figure \ref{Figure-3}(a)
is flat except a sharp rise on the left border. The flatness is
explained by the uniform
distribution of the original $p$-values corresponding to the true null
hypotheses, whereas the sharp rise is caused by the
small $p$-values corresponding to the true alternative hypotheses.
The histogram of the median aggregated $p^*$-values in Figure \ref
{Figure-3}(c) shows a
sharp rise at the
left end and has a shape symmetric about $0.5$. The approximate
symmetry arises from
the limit distribution of $p^*$-values corresponding to the true null
hypotheses [see (\ref{c.2})], whereas the
sharp rise is formed by small $p^*$-values
corresponding to the true alternative hypotheses.
Figures \ref{Figure-3}(b), (d) and (d$'$) manifest that
the $\FDR$ procedure diminishes the effectiveness in detecting the
significant regions than
the $\FDR_L$ procedure, demonstrating that the $\FDR_L$ procedure
more effectively
increases the true positive rates.
As a comparison, Figures \ref{Figure-3}(e), (f) and (f$'$) correspond
to using the mean (other than median) filter for aggregating $p$-values.
It is seen that the detections by the median and mean filters are very similar;
but compared with the mean, the median better preserves
the edge of the larger black rectangle between significant and
nonsignificant areas.
This effect gets more pronounced when $\alpha$ increases,
lending support to the ``edge preservation property'' of the median.

%
%

To evaluate the performance of Method I and Method II in estimating
$\widetilde G^*(t)$,
the bottom panels of Figure \ref{Figure-3} display the plots of
$\widehat G^*(t)$ versus $\widetilde G^*(t)$ and $\widehat G^*_c(t)$
versus $\widetilde
G^*(t)$. The agreement with 45 degree lines
well supports both estimation methods.

\begin{figure}

\includegraphics{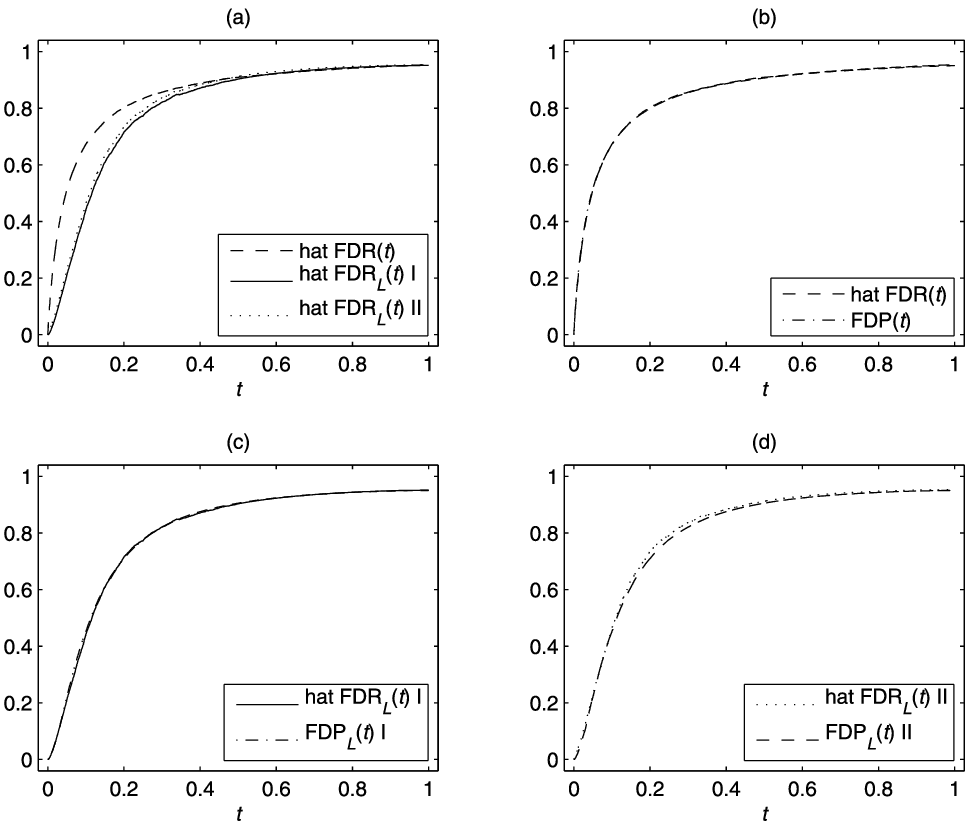}

\caption{Panel \textup{(a)}:\vspace*{1pt} compare the average values of $\widehat\FDR(t)$
and those of $\widehat\FDR_L(t)$ using Methods $\mathrm{I}$ and $\mathrm{II}$.
Panel \textup{(b)}: compare the average values of $\widehat\FDR(t)$ and those
of $\FDP(t)$.
Panel \textup{(c)}: compare the average values of $\widehat\FDR_L(t)$ using
Method $\mathrm{I}$ and those of $\FDP_L(t)$.
Panel \textup{(d)}: compare the average values of $\widehat\FDR_L(t)$ using
Method $\mathrm{II}$ and those of $\FDP_L(t)$.
Here $\lambda=0.1$.}
\label{Figure-4}
\end{figure}

To examine the overall performance of the estimated $\FDR(t)$ and
$\FDR_L(t)$ for a same threshold $t\in[0, 1]$,
we replicate the simulation $100$ times.
For notational convenience, denote by $\FDP(t)=V(t)/\{R(t)\vee1\}$ and
$\FDP_L(t)=V^*(t)/\{R^*(t)\vee1\}$ the false discovery proportions of the
$\FDR$ and $\FDR_L$ procedures, respectively.
The average values (over $100$ data) of $\widehat\FDR(t)$ and $\widehat
\FDR
_L(t)$ at each point
$t$ are plotted in Figure \ref{Figure-4}(a).
%
%
It is clearly observed that $\widehat\FDR_L(t)$ using both Methods
$\mathrm{I}$ and $\mathrm{II}$ is below $\widehat\FDR(t)$,
demonstrating that
the $\FDR_L$ procedure produces the estimated\vspace*{1pt} false discovery rates
lower than those of the $\FDR$ procedure.
Meanwhile, Figure \ref{Figure-4} compares the average values of
$\FDP(t)$ and those of $\widehat\FDR(t)$ in panel (b), and the average
values of
$\FDP_L(t)$ using Methods $\mathrm{I}$ and $\mathrm{II}$ and those
of $\widehat\FDR_L(t)$ in panels (c) and (d), respectively.
For each procedure, the two types of estimates are very close to each
other, lending
support to the estimation procedure in Section~\ref{sec-3.4}.

\subsubsection{Sensitivity and specificity}

To further study the relative performance of the $\FDR$ and $\FDR_L$
procedures,
we adopt two widely used performance measures,
\begin{eqnarray*}
\mathrm{sensitivity}
&\equiv&
\cases{
{S(t_\alpha(\widehat\FDR))}/{n_1}, &\quad \mbox{for the $\FDR$ procedure},
\cr
{S^*(t_\alpha(\widehat\FDR_L))}/{n_1}, &\quad \mbox{for the $\FDR_L$
procedure},}
\\
\mathrm{specificity}
&\equiv&
\cases{
{U(t_\alpha(\widehat\FDR))}/{n_0}, &\quad \mbox{for the $\FDR$ procedure},
\cr
{U^*(t_\alpha(\widehat\FDR_L))}/{n_0}, &\quad \mbox{for the $\FDR_L$
procedure},}
\end{eqnarray*}
for summarizing the discriminatory power of
a diagnosis procedure,
where
$S(t) = \sum_{i=1}^n \ID\{H_0(i)$ is false, and $p_i \le t\}$,
$U(t) = \sum_{i=1}^n \ID\{H_0(i)$ is true, and $p_i > t\}$,
$S^*(t) = \sum_{i=1}^n \ID\{H_0(i)$ is false, and $p_i^* \le t\}$
and
$U^*(t) = \sum_{i=1}^n \ID\{H_0(i)$ is true, and $p_i^* > t\}$.
Here, the sensitivity and specificity measure the strengths for correctly
identifying the alternative and the null hypotheses, respectively.

%
%
\begin{figure}

\includegraphics{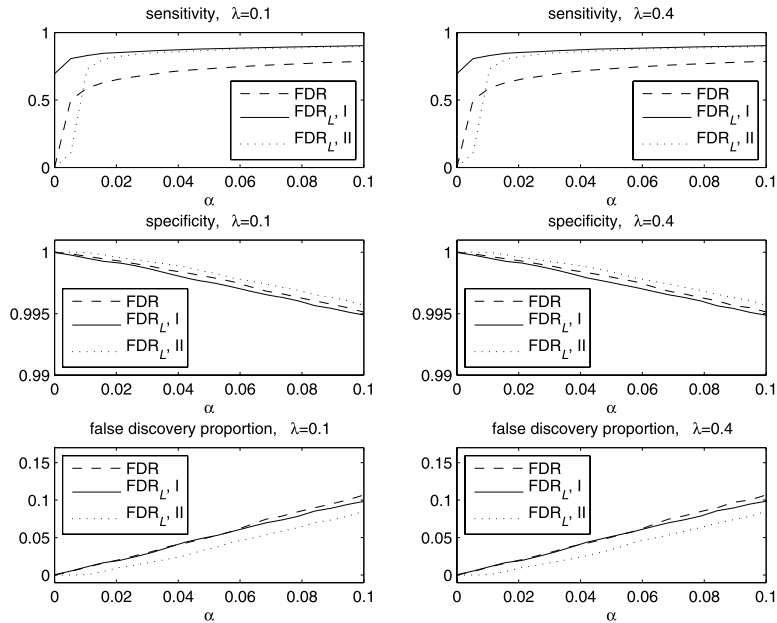}

\caption{Comparison of the average sensitivity (top
panels), average specificity (middle panels) and average false
discovery proportion (bottom panels). Left panels: $\lambda=
0.1$. Right panels: $\lambda= 0.4$.}\label{Figure-5}
\end{figure}

Following Section \ref{sec-5.1}, we randomly generate $100$ sets of
simulated data and perform $\FDR$ and $\FDR_L$ procedures for each dataset,
with the control levels $\alpha$ varying from $0$ to $0.1$.
The left panel of Figure \ref{Figure-5} corresponds to $\lambda=0.1$,
whereas the right panel corresponds to $\lambda=0.4$.
In either case, we observe that the average sensitivity (over the
datasets) of the $\FDR_L$ procedure using Method I
is consistently higher than that of the $\FDR$ procedure,
whereas the average specificities of both procedures approach one and
are nearly indistinguishable.
In addition, the bottom panels indicate that the $\FDR$ procedure
yields larger (average) false discovery proportions
than the $\FDR_L$ procedure.
It is apparent that the results in Figure~\ref{Figure-5} are not very
sensitive to the choice of $\lambda$.
Unless otherwise stated, $\lambda=0.1$ will be used throughout the
rest of the numerical work.

\subsection{Example $2$: More strongly correlated case} \label{sec-5.2}
We consider a dataset generated according to the same
model (\ref{e.1}) as in Example 1, but with more strongly correlated
errors,\vspace*{2pt}
by taking
$\varepsilon(i,j)=\sum_{i=0}^{6} \sum_{j=0}^{6} e(i,j)/7$,
where $\{e(i,j)\}_{i,j=0}^{264}$ are i.i.d. $N(0,1)$.
As seen from the figure in \citet{ZFY10},
both $\FDR$ and $\FDR_L$ (using Methods I and II) procedures perform
worse with strongly-correlated data
than with low-correlated data (given in Figure \ref{Figure-3}). However,
there are no adverse effects by applying $\FDR_L$ to more strongly
correlated data,
and
Method I continues to be comparable with Method II for the $\FDR_L$ procedure.

\subsection{Example $3$: Large proportion of boundary grid points}
\label{sec-5.3}
The efficacy of the $\FDR_L$ procedure is illustrated in the figure of
\citet{ZFY10}
by a simulated dataset
generated according to the same model (\ref{e.1}) as in Example 1, but with
a large proportion of boundary grid points, where $\mu(i,j) = 0$ for
$(i,j)\in\calV_0$ and
$\mu(i,j) = 4$ for $(i,j)\in\calV_1$. Similar plots using $\mu(i,j)
= 2$ for $(i,j)\in\calV_1$ are obtained and thus omitted.
Again, there is no adverse effect of using $\FDR_L$ to detect dense or
weak signals.

\section{Simulation study: 3D dependent data} \label{sec-6}

We apply the $\FDR$ and $\FDR_L$ procedures to detect activated brain regions
of a simulated brain fMRI dataset, which is both spatially and
temporally correlated.
The experiment
design, timings and size are exactly the same as those of the real fMRI
dataset in
Section \ref{sec-7}.
%
%
\begin{figure}[b]
\begin{tabular}{@{}cc@{}}

\includegraphics{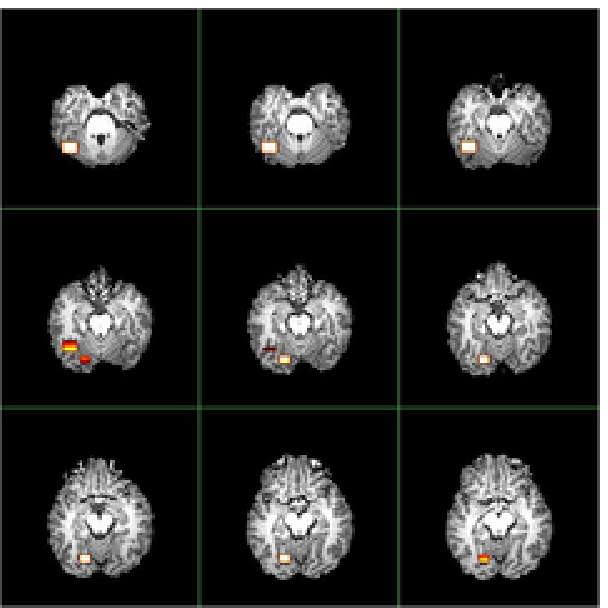}
 & \includegraphics{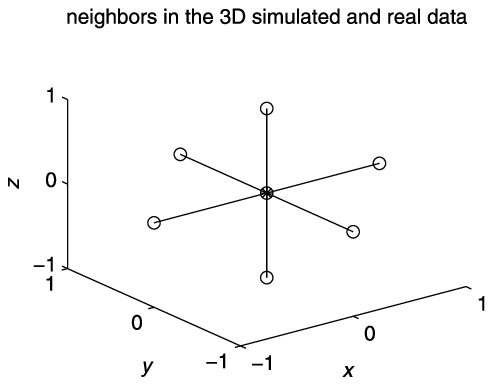}
\end{tabular}
\caption{Left panel: true activated brain regions (denoted
by hot color) for the simulated $\fMRI$ dataset. Right panel:
neighbors of a point at $(x,y,z)$ used in the $\FDR_L$ procedure for
$\mathrm{3D}$ simulated and real data.} \label{Figure-8}
\end{figure}
The data are generated from a semi-parametric model similar to that in
Section 5.2 of \citet{ZY08}. (They demonstrated that the
semi-parametric model gains more flexibilities than existing parametric
models.) The left panel of Figure \ref{Figure-8} contains $9$ slices
(corresponding to the 2D axial view) which highlight two activated
brain regions involving $91$ activated brain voxels. The neighborhood
used in the $\FDR_L$ procedure is illustrated in the right panel of
Figure \ref{Figure-8}.

%
%
\begin{figure}

\includegraphics{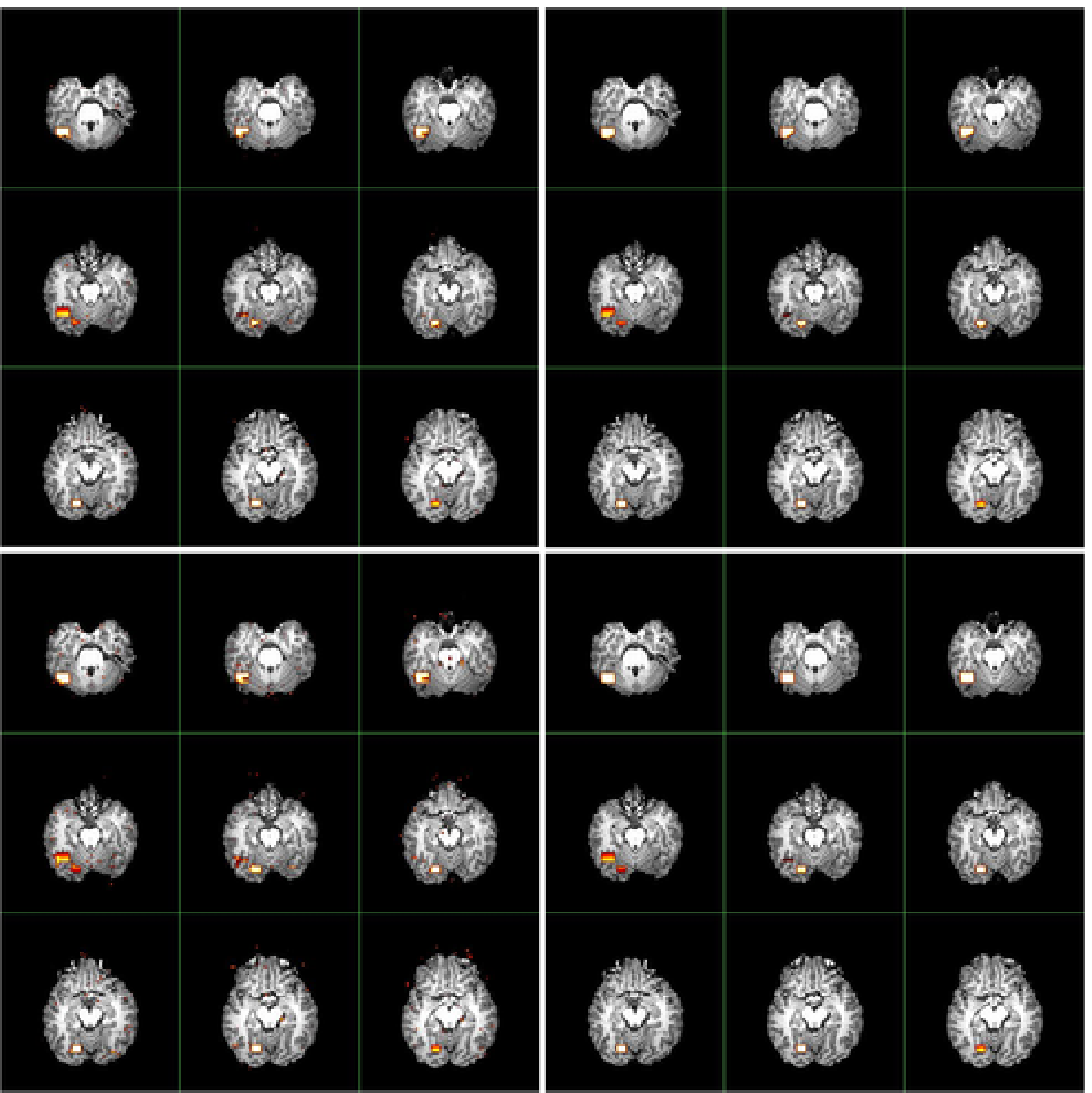}

\caption{Comparison of activated brain regions detected for
the simulated $\fMRI$ dataset using the conventional $\FDR$ approach
(on the left) and the proposed $\FDR_L$ procedure (on the
right) using Method $\mathrm{I}$. Top panels: $\chiK$. Bottom
panels: $\chiK_\bc$. Here $\alpha=0.05$.} \label{Figure-9}
\end{figure}

Figure \ref{Figure-9} compares the activated brain regions identified by
the $\FDR$ (in the left panels) and $\FDR_L$ (in the right panels) procedures.
Owing to the wealth of data, and for purposes of computational simplicity,
results using Method I of $\FDR_L$ are presented.
Voxel-wise\vspace*{1pt} inactivity is tested with
the semi-parametric test statistics
$\chiK=(A \widehat\bh)^T \{A(\widetilde\bS^T \widehat R^{-1} \widetilde
\bS
)^{-1} A^T \}^{-1} (A \widehat\bh)/\break
\{\widehat\br^T \widehat R^{-1} \widehat\br/(n-rm)\}$ (in\vspace*{1pt} the top panels)
and $\chiK_\bc=(A \widehat\bh_\bc)^T \{A(\widetilde\bS^T \widehat R^{-1}
\widetilde
\bS)^{-1}\* A^T \}^{-1} (A \widehat\bh_\bc)/
\{\widehat\br_\bc^T \widehat R^{-1} \widehat\br_\bc/(n-rm)\}$ (in the
bottom panels)
whose notation was given and asymptotic $\chi^2$ distributions were
derived in \citet{ZY08}.
The control level is $0.05$.
Inspection of Figure~\ref{Figure-9} reveals that
$\chiK$ and $\chiK_\bc$ locate both active regions.
In particular,
using the $\FDR$ procedure, both methods detect more than $200$ voxels
(which are visible when zooming the images), many of which are
falsely discovered.
When applying the $\FDR_L$ procedure,
$\chiK$ detects $82$ voxels, whereas
$\chiK_\bc$ detects $90$ voxels.
Thus the $\FDR_L$ procedure
reduces the number of tiny scattered false findings,
gaining more accurate detections than the $\FDR$ procedure.

%
%
\begin{figure}

\includegraphics{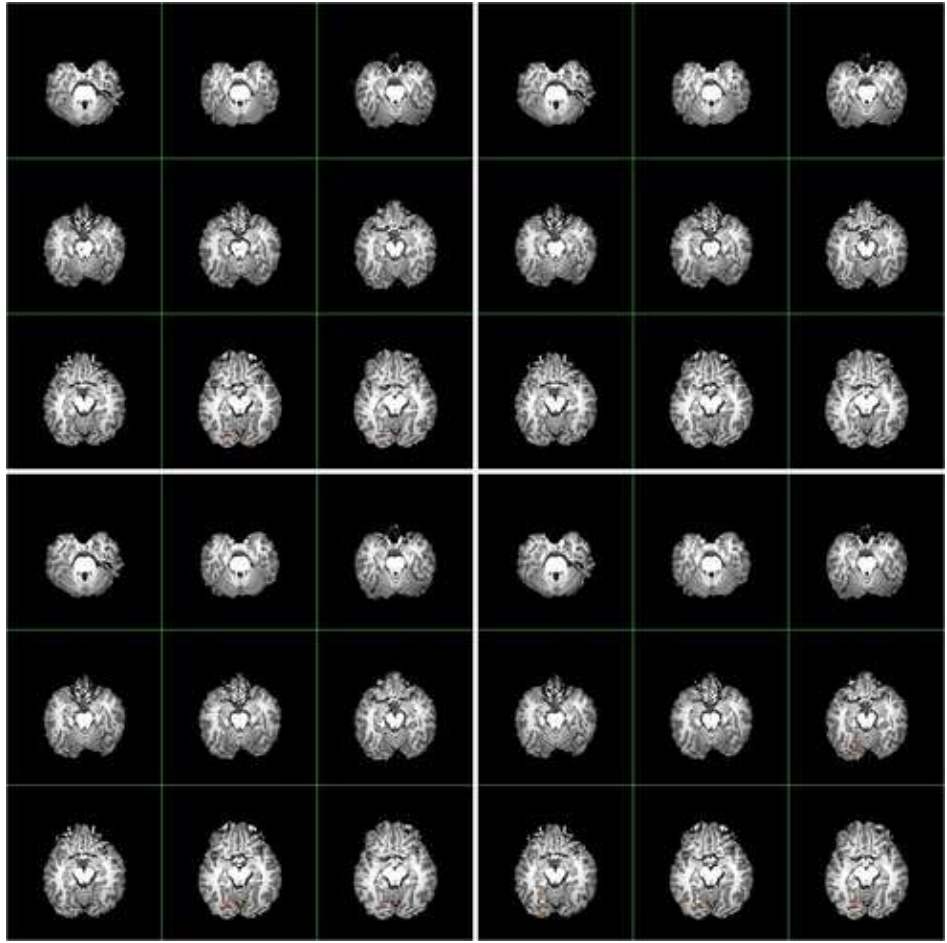}

\caption{Comparison of activated brain regions detected for
the simulated $\fMRI$ dataset using the conventional $\FDR$ approach
(on the left) and the proposed $\FDR_L$ procedure (on the
right) using Method $\mathrm{I}$. Top panels: $\AFNI$. Bottom
panels: $\FSL$. Here $\alpha=0.05$.}\label{Figure-10}
\end{figure}

As a comparison, the detection results by popular software AFNI [\citet{C96}]
and FSL [\citet{Setal04} and
\citet{Wetal01}] are given in Figure \ref{Figure-10}.
We observe that both AFNI and FSL
fail to locate one activated brain area, and that the other region,
though correctly detected, has appreciably reduced size relative to
the actual size. This detection bias is due to the stringent
assumptions underlying AFNI and FSL in modeling fMRI data:
the Hemodynamic Response Function (HRF) in FSL is specified as the
difference of two gamma functions, and the
drift term in AFNI is specified as a quadratic polynomial. As anticipated,
applying the $F$ distributions restricted to parametric models to
specify the
distributions of test statistics in AFNI and FSL
leads to bias, which in turn gives
biased calculations of $p$-values and $p^*$-values.
In this case, the detection performances of both the $\FDR$ and $\FDR
_L$ procedures
deteriorate, and the $\FDR_L$ procedure does not improve the
performance of the $\FDR$ procedure.
See Table \ref{Table-3} for a more detailed comparison.

To reduce modeling bias, for
applications to the real fMRI dataset in Section~\ref{sec-7}, we
will only employ the semi-parametric test statistics $\chiK$ and
$\chiK_\bc$.
It is also worth distinguishing between the computational aspects
associated with the $\FDR_L$ procedure:
this paper uses (\ref{c.3}) for the null distribution of $p^*$-values,
whereas \citet{ZY08} used the normal approximation approach in
Section \ref{sec-3.2}.

\section{Functional neuroimaging example} \label{sec-7}

In an emotional control study, subjects saw a series of negative or
positive emotional images, and
were asked to either suppress or enhance their emotional responses to
the image, or to simply attend
to the image.
The sequence of trials was randomized.
The time between successive trials also varied.
The size of the whole brain dataset is $64\times64\times30$. At
each voxel, the time series has $6$ runs, each containing $185$
observations with a time resolution of $2$ seconds.
For details of the dataset, please refer to \citet{ZY08}.
The study aims to estimate the BOLD (Blood Oxygenation
Level-Dependent) response to each of the trial types for $1$--$18$
seconds following the image onset. We analyze
the fMRI dataset containing one subject. The length of the estimated
HRF is set equal to $18$.
Again, the neighborhood used in the $\FDR_L$ procedure is illustrated in
the right panel of Figure \ref{Figure-8}.

%
%
\begin{table}
\tabcolsep=0pt
\caption{Comparing ${\FDR}$ and ${\FDR_L}$ procedures}
\label{Table-3}
\begin{tabular*}{\tablewidth}{@{\extracolsep{\fill}}lld{3.4}d{3.4}d{2.4}d{2.4}@{}}
\hline
& & \multicolumn{4}{c@{}}{\textbf{Test methods}} \\
[-4pt]
& & \multicolumn{4}{c@{}}{\hrulefill}\\
& \multicolumn{1}{c}{\textbf{Multiple comparison}}
& \multicolumn{1}{c}{$\bolds{\chiK}$} & \multicolumn{1}{c}{$\bolds{\chiK_\bc}$}
& \multicolumn{1}{c}{\textbf{AFNI}} & \multicolumn{1}{c@{}}{\textbf{FSL}} \\
\hline
\# of detected voxels & $\FDR$ & 276 & 870 & 16 & 6 \\
& $\FDR_L$, Method I & 82 & 90 & 2 & 11 \\ [2pt]
False discovery proportion & $\FDR$ & 0.6993 & 0.9000 & 0.5625 & 0 \\
& $\FDR_L$, Method I & 0 & 0 & 0.5000 & 0 \\ [2pt]
Sensitivity & $\FDR$ & 0.9121 & 0.9560 & 0.0769 & 0.0659 \\
& $\FDR_L$, Method I & 0.9011 & 0.9890 & 0.0110 & 0.1209 \\ [2pt]
Specificity & $\FDR$ & 0.9921 & 0.9678 & 0.9996 & 1.0000 \\
& $\FDR_L$, Method I & 1.0000 & 1.0000 & 0.9997 & 1.0000 \\
\hline
\end{tabular*}
\end{table}

%
%
\begin{figure}

\includegraphics{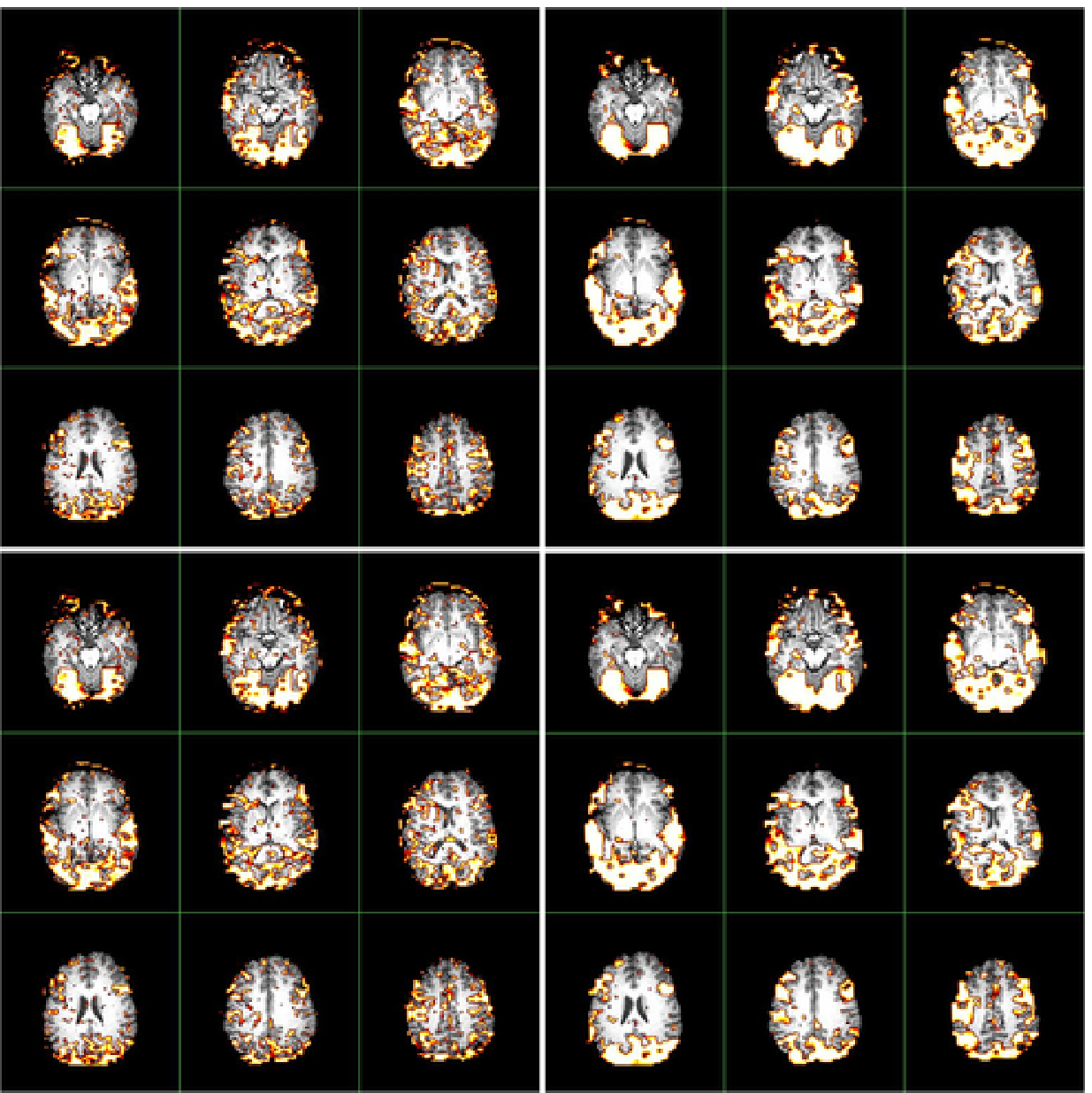}

\caption{Comparison of activated brain regions detected for
the real $\fMRI$ dataset using the conventional $\FDR$ approach (on
the left) and the proposed $\FDR_L$ procedure (on the right)
using Method~$\mathrm{I}$. Top panels: $\chiK$. Bottom panels:
$\chiK_\bc$. Here $\alpha=0.001$.} \label{Figure-11}
\end{figure}

A comparison of the activated brain regions
using the $\FDR$ and $\FDR_L$ procedures
is visualized in Figure
\ref{Figure-11}.
The level $0.001$ is used to carry out the multiple comparisons.
The conventional $\FDR$ procedure finds more tiny scattered active
voxels, which are more likely to be
falsely discovered. In contrast,
the $\FDR_L$ procedure finds activation in much more clustered
regions of the brain.

\section{Discussion} \label{sec-8}

This paper proposes the $\FDR_L$ procedure
to embed the structural spatial information of $p$-values
into the conventional $\FDR$ procedure for large-scale imaging data
with a spatial structure.
This procedure provides the standard $\FDR$ procedure with the ability
to perform better on spatially aggregated $p$-values.
Method I and Method II have been developed for making statistical inference
of the aggregated $p$-values under the null.
Method I gains remarkable computational superiority, particularly for large/huge
imaging datasets, when the $p^*$-values under the null are not too skewed.
Furthermore, we provide a better understanding of
a ``\textit{lack of identification phenomenon}'' ($\LIP$) occurring in the
$\FDR$ procedure. This study indicates that the $\FDR_L$ procedure
alleviates the extent of the problem and
can adopt control levels
much smaller than those of the $\FDR$ procedure without excessively
encountering the $\LIP$, thus substantially facilitating
the selection of more stringent
control levels.

As discussed in \citet{O05} and \citet{LS08}, a key issue
with the
dependencies between the hypotheses tests is the inflation of the
variance of significance measures in $\FDR$-related work.
Indeed, similar to $\FDR$, the $\FDR_L$ procedure (using Methods I
and II) performs less well with highly-correlated data
than with the low-correlated data. Detailed investigation of the
variance of $\FDR_L$ will be given in future study.

Other ways of exploring spatially neighboring information are certainly possible
in multiple comparison. For example,
the median operation applied to $p$-values can be replaced by
the averaging, kernel smoothing, ``majority vote'' and
edge preserving smoothing techniques [Chu et al. (\citeyear{Cetal98})].
Hence, taking the median is not the unique way to aggregate $p$-values.
On the other hand, compared with the mean,
the median is more robust, computationally simpler and does not depend
excessively on
the spatial co-ordinates, especially on the boundaries between
significant and nonsignificant
regions, as observed in Figures \ref{Figure-3}(d) and (f).
An exhaustive comparison is beyond the scope of the current paper
and we leave this for future research.

\begin{appendix}\label{app}
\section{\texorpdfstring{Proofs of Theorems 4.1--4.3}{Proofs of Theorems
4.1--4.3}}
\label{appA}

We first impose some technical assumptions, which are not the weakest possible.
Detailed proofs of Theorems \ref{Thm-1}--\ref{Thm-3} are given in
\citet{ZFY10}.

{\renewcommand{\theCondition}{A}
\begin{Condition}\label{CondA}
\smallskipamount=0pt
\begin{longlist}
\item[A0.]\hypertarget{CondA0}
The neighborhood size $k$ is an integer not depending on $n$.

\item[A1.]\hypertarget{CondA1} $\lim_{n\to\infty} {n_0}/{n} = \pi_0$ exists and $\pi
_0 < 1$.

\item[A2.]\hypertarget{CondA2}
$\lim_{n \to\infty} {V^*(t)}/{n_0} = G_0^*(t)$ and
$\lim_{n \to\infty} {S^*(t)}/{n_1} = G_1^*(t)$ almost surely for
each $t\in(0,1]$, where $G_0^*$ and $G_1^*$
are continuous functions.

\item[A3.]\hypertarget{CondA3} $0<G_0^*(t)\le G^{*\infty}(t)$ for each
$t\in(0,1]$.\vspace*{2pt}

\item[A4.]\hypertarget{CondA4} $\sup_{t\in(0,1]} |\widehat G^*(t)-G^{*\infty}(t)| = o(1)$ almost
surely as $n\to\infty$.
\end{longlist}
\end{Condition}}

\vspace*{-10pt}
\section{\texorpdfstring{Proofs of Theorems 4.4 and 4.5}{Proofs of Theorems 4.4 and 4.5}}
\label{appB}

\subsection{\texorpdfstring{Proof of Theorem \protect\ref{Thm-4}}{Proof of Theorem 4.4}}
By the assumptions and $H_1(v)$, we see that the $p$-value has the expression,
$p(v)
= 1-F_0(T(v))$.
Thus,
the distribution function of $p(v)$ corresponding to the true $H_0(v)$ is
$G_0(t)=t$ for $0<t<1$ and (\ref{b.3}) gives
$\widehat{\FDR}^{\infty}(t)=\frac{\pi_0+\pi_1\{1-G_1(\lambda)\}
/(1-\lambda)} {\pi_0+\pi_1 G_1(t)/t}$.
Also,
the distribution function of $p(v)$ corresponding to the true $H_1(v)$
is given by
%
%
\begin{equation} \label{B.1}
G_1(t)
= 1-F_1\bigl(F_0^{-1}(1-t)\bigr).\vadjust{\goodbreak}
\end{equation}
Likewise, using
(\ref{c.2}), it follows that with probability one,
%
%
\begin{eqnarray} \label{B.2}
G_0^*(t)
&=& \lim_{n\to\infty}\frac{V^*(t)}{n_0}\nonumber\\
&=&\lim_{n\to\infty}\frac{\sum_{v\in\calV_0^{(0)}}\ID\{p^*(v)\le
t\}}{\# \calV_0^{(0)}} \cdot\lim_{n\to
\infty}\frac{\# \calV_0^{(0)}}{n_0} \nonumber\\
&&{} + \lim_{n\to\infty}\frac{\sum_{v\in\calV_0 \setminus\calV
_0^{(0)}}\ID\{p^*(v) \le
t\}}{n_0} \\
&=& P\{p^*(v) \le t\}\qquad \mbox{with } v \in\calV_0^{(0)} \nonumber\\
&=& G^{*\infty}(t) = B_{(k+1)/2, (k+1)/2}(t),\nonumber
\end{eqnarray}
the cumulative distribution function of a ${\Betadist
}((k+1)/2,(k+1)/2)$ random variable and
%
%
\begin{equation} \label{B.3}
G_1^*(t) =\lim_{n\to\infty} {S^*(t)}/{n_1}= B_{(k+1)/2,
(k+1)/2}(G_1(t)).
\end{equation}
Applying (\ref{B.2}) and (\ref{d.1}) gives $\widehat\FDR_L^{\infty}(t) =
\frac{\pi_0+\pi_1\{1-G_1^*(\lambda)\}/\{1- G_0^*(\lambda)\}}{\pi
_0+\pi_1 G_1^*(t)/G_0^*(t)}$.

{Part I.} For the $\FDR$ procedure, note that $\widehat{\FDR}^{\infty}(t)$
is a decreasing function of $G_1(t)/t$. Applying L'Hospital's rule and
the fact $\lim_{t\to0+} G_1(t) =
0$,
%
%
\begin{equation} \label{B.4}
\lim_{t\to0+} \frac{G_1(t)}{t}
=
\lim_{t\to0+}\frac{f_1(F_0^{-1}(1-t))}{f_0(F_0^{-1}(1-t))}
= \lim_{x\to x_0-}\frac{f_1(x)}{f_0(x)} = \infty,
\end{equation}
where $x = F_0^{-1}(1-t)$. Thus, $\sup_{0 < t \le1}G_1(t)/t = \infty$,
which together with $\widehat{\FDR}^{\infty}(t)$ shows $\alpha_{\infty
}^{\FDR} =0$ for the $\FDR$
procedure.

For the $\FDR_L$ procedure,
applying (\ref{B.2}) and (\ref{B.3}), we get
%
%
\begin{eqnarray}\qquad
\label{B.5}
\frac{dG_0^*(t)}{dt}
&=& \frac{dG^{*\infty}(t)}{dt} = \frac{k!}{[\{(k-1)/2\}
!]^2}t^{(k-1)/2}(1-t)^{(k-1)/2},\\
\label{B.6}
\frac{dG_1^*(t)}{dt}
&=& \frac{k!}{[\{(k-1)/2\}!]^2}G_1(t)^{(k-1)/2}\{1-G_1(t)\}^{(k-1)/2}\,
\frac{d G_1(t)}{dt}.
\end{eqnarray}
Note that $\widehat\FDR_L^{\infty}(t)$
is a decreasing function of $G_1^*(t)/G_0^*(t)$. Since\break $\lim_{t\to0+}
G_1^*(t)=0$
and $\lim_{t\to0+} G_0^*(t)=0$,
%
%
\begin{eqnarray} \label{B.7}
\lim_{t\to0+} \frac{G_1^*(t)}{G_0^*(t)}
&=& \lim_{t\to0+}
\frac{{dG_1^*(t)}/{dt}}{{dG_0^*(t)}/{dt}}\nonumber\\[-8pt]\\[-8pt]
&=& \lim_{t\to0+}
\biggl\{
\frac{G_1(t)}{t}\cdot
\frac{1-G_1(t)}{1-t}\biggr\}^{(k-1)/2}\,
\frac{dG_1(t)}{dt},\nonumber
\end{eqnarray}
which together with (\ref{B.4})
shows
$
\lim_{t\to0+} G_1^*(t)/G_0^*(t) = \infty$.
Thus,
\[
\sup_{0 < t \le1} G_1^*(t)/G_0^*(t)=\infty,
\]
that is,
$\alpha_{\infty}^{\FDR_L} = 0$ for the $\FDR_L$ procedure.

{Part II.}
Following $\widehat{\FDR}^{\infty}(t)$ and $\widehat\FDR_L^{\infty
}(t)$, we
immediately conclude that
$\alpha_{\infty}^{\FDR} \ne0$ if
%
%
\begin{equation} \label{B.8}
\sup_{0 < t \le1}G_1(t)/t<\infty,
\end{equation}
and that $\alpha_{\infty}^{\FDR_L} \ne0$ if
%
%
\begin{equation}\label{B.9}
\sup_{0 < t \le1}G_1^*(t)/G_0^*(t)<\infty.
\end{equation}

We first verify (\ref{B.8}) for the $\FDR$ procedure. Assume (\ref
{B.8}) fails, that is, $\sup_{0<t \le1}G_1(t)/t=\infty$. Note that
for any $\delta>0$, the function
$G_1(t)/t$, for $t\in[\delta, 1]$, is continuous and bounded away from
$\infty$, thus, $\sup_{0<t \le1}G_1(t)/t=\infty$ only if there exists
a sequence $t_1>t_2>\cdots>0$, such that $\lim_{m\to\infty} t_m =
0$ and $\lim_{m\to\infty} G_1(t_m)/t_m = \infty$. For each $m$,
recall that
both $G_1(t)$ and $t$ are continuous on $[0,t_m]$, and differentiable
on $(0,t_m)$. Applying
Cauchy's mean-value theorem, there exists $\xi_m \in(0,t_m)$ such that
$
{G_1(t_m)}/{t_m}
=
\{G_1(t_m)-G_1(0)\}/{(t_m-0)}
=
\frac{d G_1(t)}{dt}|_{t=\xi_m}.
$
Since $\lim_{m\to\infty} G_1(t_m)/t_m = \infty$, it follows that
%
%
\begin{equation} \label{B.10}
\limsup_{t\to0+}\frac{d
G_1(t)}{dt} = \infty.
\end{equation}
On the other hand, the condition $\limsup_{x\to x_0-}\frac
{f_1(x)}{f_0(x)} < \infty$ indicates that
%
%
\begin{equation} \label{B.11}\quad
\limsup_{t\to0+} \frac{d G_1(t)}{dt}
= \limsup_{t\to0+}\frac{f_1(F_0^{-1}(1-t))}{f_0(F_0^{-1}(1-t))}
= \limsup_{x\to x_0-}\frac{f_1(x)}{f_0(x)} <\infty,
\end{equation}
where $x = F_0^{-1}(1-t)$. Clearly, (\ref{B.11}) contradicts (\ref{B.10}).

Next, we show (\ref{B.9}) for the $\FDR_L$ procedure. Combining (\ref
{B.7}), (\ref{B.8}) and
(\ref{B.11}), the result follows. This completes the proof.

\subsection{\texorpdfstring{Proof of Theorem \protect\ref{Thm-5}}{Proof of Theorem 4.5}}
We first show Lemma \ref{Lemma-1}.
\begin{lemma} \label{Lemma-1}
Let $B(t)$ be the cumulative distribution function of a ${\Betadist
}(a,a)$ random variable, where $a>1$ is a real number. Then
$\mathrm{I}$
for $t\in(0, 0.5)$, $B(t)/t$ is a strictly increasing function and $B(t)<t$;
$\mathrm{II}$
for $t\in(0.5, 1)$, $B(t)>t$;
$\mathrm{III}$
for $t_1 \in(0, 0.5]$ and $t_2\in[t_1, 1]$, $B(t_1)/t_1\le B(t_2)/t_2$.
\end{lemma}
\begin{pf}
Let $\Gamma(\cdot)$ denote the Gamma function. It is easy to see that
%
%
\begin{equation} \label{B.12}
B''(t)= {\Gamma(2a)}/{\{\Gamma(a)\}^2}(a-1)t^{a-2}(1-t)^{a-2}(1-2t).
\end{equation}

To show part I, define $F_1(t)=B(t)/t$.
Then $F_1'(t)=\{B'(t)t-B(t)\}/t^2$, where
$\frac{d\{B'(t)t- B(t)\}}{dt} = B''(t) t$.
For $t\in(0, 0.5)$, (\ref{B.12}) indicates
$B''(t)>0$, that is, $B'(t)t- B(t)$ is strictly increasing,
implying $B'(t)t-B(t)> B'(0)0-B(0)=0$.
Hence for $t\in(0, 0.5)$, $B(t)/t$ is strictly increasing,
and therefore $B(t)/t < B(0.5)/0.5 = 1$.

For part II, define $F_2(t)=B(t)-t$. Then
$F_2''(t)=B''(t)$. By (\ref{B.12}), $B''(t)<0$ for $t\in(0.5, 1)$,
thus $F_2(t)$ is strictly concave, giving
$F_2(t)>\max\{F_2(0.5)$, $F_2(1)\}=0$.

Last, we show part III.
For $t_2 \in[t_1, 0.5]$, part I indicates that
$B(t_1)/t_1\le B(t_2)/t_2$;
for $t_2 \in[0.5, 1]$, part II indicates that
$B(t_2)/t_2\ge1$ which,
combined with $B(t_1)/t_1\le1$ from part I,
yields $B(t_1)/t_1 \le B(t_2)/t_2$.
\end{pf}

We now prove Theorem \ref{Thm-5}. It suffices to show that
%
%
\begin{eqnarray}
\label{B.13}
\{1-G_1(\lambda)\}/{(1-\lambda)}
&\ge& \{1-G_1^*(\lambda)\}/\{1-G_0^*(\lambda)\},\\
\label{B.14}
\sup_{0< t\le1} {G_1(t)}/{t}
&\le& \sup_{0<t \le1} {G_1^*(t)}/{G_0^*(t)}.
\end{eqnarray}
Following (\ref{B.5}) and (\ref{B.6}), for $0\le t\le1$,
%
%
\begin{equation} \label{B.15}
G_1^*(t)
=G_0^*(G_1(t)).
\end{equation}
Applying
(\ref{B.15}),
(\ref{B.1}),
$1-F_0(F_1^{-1}(0.5))\le\lambda$
and part II of Lemma \ref{Lemma-1} yields $G_1(\lambda) \le
G_1^*(\lambda)$;
applying $\lambda\le0.5$ and part
I of Lemma \ref{Lemma-1} implies $\lambda\ge G_0^*(\lambda)$. This
shows (\ref{B.13}).

To verify (\ref{B.14}), let $M=\sup_{0< t\le1} {G_1(t)}/{t}$.
Since $G_1(1)/1 = 1$, we have $M\ge1$ which will be discussed in two cases.
Case 1: if $M=1$, then
%
%
\begin{equation} \label{B.16}
\sup_{0<t \le1}
\frac{G_1^*(t)}{G_0^*(t)} \ge\frac{G_1^*(1)}{G_0^*(1)} = 1 =
\sup_{0< t\le1} \frac{G_1(t)}{t}.
\end{equation}

Case 2: if $M>1$, then there exists $t_0\in[0,1]$ and $t_n \in(0,1)$
such that $\lim_{n\to\infty}t_n=t_0$, and
%
%
\begin{equation} \label{B.17}
\lim_{n\to\infty} G_1(t_n)/t_n =\sup_{0<t\le1}
G_1(t)/t = M>1.
\end{equation}
Thus, there exists $N_1$ such that for all $n>N_1$,
%
%
\begin{equation} \label{B.18}
G_1(t_n)>t_n.
\end{equation}
Cases of $t_0=1$, $t_0=0$ and $t_0\in(0,1)$ will be discussed separately.
First, if \mbox{$t_0=1$}, then $M=\lim_{n\to\infty} G_1(t_n)/t_n =
\lim_{n\to\infty}G_1(t_n) \le1$, which contradicts
(\ref{B.17}). Thus $t_0<1$.
Second, if $t_0=0$, then there exists $N_2$ such that $t_n<0.5$ for all
$n>N_2$. Thus
for all $n>N\equiv\max\{N_1,N_2\}$, applying
(\ref{B.15}),
(\ref{B.18}) and
part III of Lemma \ref{Lemma-1}, we have that
\[
\frac{G_1^*(t_n)}{G_1(t_n)} = \frac{G_0^*(G_1(t_n))}{G_1(t_n)}\ge
\frac{G_0^*(t_n)}{t_n}.
\]
This together with (\ref{B.17}) shows
%
%
\begin{equation} \label{B.19}\qquad
\sup_{0<t\le1}\frac{G_1^*(t)}{G_0^*(t)}\ge\limsup_{n\to\infty}
\frac{G_1^*(t_n)}{G_0^*(t_n)} \ge
\lim_{n\to\infty}\frac{G_1(t_n)}{t_n} = M =\sup_{0<t\le
1}\frac{G_1(t)}{t}.
\end{equation}
Third, for $t_0 \in(0,1)$, since both $F_0$ and $F_1$ are
differentiable and $f_0$ is supported in a single interval,
$G_1(t)/t=\{1-F_1(F_0^{-1}(1-t))\}/t$ is differentiable in $(0,1)$.
Thus,
%
%
\begin{equation} \label{B.20}
\sup_{0<t\le1} {G_1(t)}/{t} = {G_1(t_0)} / {t_0} = M
\end{equation}
and $\frac{d\{G_1(t)/t\}}{dt}|_{t=t_0} = 0$.
Notice
%
%
\begin{eqnarray} \label{B.21}
\frac{d\{G_1(t)/t\}}{dt}\bigg|_{t=t_0}&=&
\frac{({dG_1(t)}/{dt})|_{t=t_0}-G_1(t_0)/t_0}{t_0} \nonumber\\[-8pt]\\[-8pt]
&=& \frac{
({dG_1(t)}/{dt})|_{t=t_0}-M}{t_0} =
0.\nonumber
\end{eqnarray}
If $t_0>0.5$, then $F_0^{-1}(1-t_0)\le F_0^{-1}(0.5)$.
By (\ref{B.4}) and the assumption on $f_0$ and $f_1$,
$
\frac{dG_1(t)}{dt}|_{t=t_0} =
f_1(F_0^{-1}(1-t_0))/f_0(F_0^{-1}(1-t_0)) \le1$,
which contradicts (\ref{B.21}). Thus,
$
0< t_0\le0.5$.
This together with
(\ref{B.15}),
(\ref{B.20}),
and part III of Lemma \ref{Lemma-1} gives
\[
\frac{G_1^*(t_0)}{G_1(t_0)} = \frac{G_0^*(G_1(t_0))}{G_1(t_0)} \ge
\frac{G_0^*(t_0)}{t_0}.
\]
This, together with (\ref{B.20}), shows
%
%
\begin{equation} \label{B.22}
\sup_{0<t\le1}\frac{G_1^*(t)}{G_0^*(t)} \ge
\frac{G_1^*(t_0)}{G_0^*(t_0)} \ge\frac{G_1(t_0)}{t_0} = M =
\sup_{0<t\le1}\frac{G_1(t)}{t}.
\end{equation}
Combining (\ref{B.16}), (\ref{B.19}) and (\ref{B.22}) completes
the proof.

\section{$\alpha_{\infty}^{\FDR}$ and $\alpha_{\infty}^{\FDR_L}$
in Table 2 of Section 4.3}\label{appC}

Before calculating $\alpha_{\infty}^{\FDR}$ and $\alpha_{\infty
}^{\FDR_L}$, we first present two lemmas.
\begin{lemma} \label{Lemma-2}
Let $f(x)$ and $g(x)$ be differentiable functions in $x\in(a,b)
\subseteq\RR$.
Suppose that $g(x)\ne0$ for $x\in(a,b)$, and
$f(x)/g(x)$ is a nonincreasing function of $x$.
For any $C \in(0,\infty)$ such that $g(x)+C \ne0$,
if ${d f(x)}/{dx} \le{d g(x)}/{dx}$ for all $x\in(a,b)$, then
$\{f(x)+C\}/\{g(x)+C\}$ is a decreasing function in $x\in(a,b)$.
\end{lemma}
\begin{pf}
The proof is straightforward and is omitted.
\end{pf}
\begin{lemma} \label{Lemma-3}
The function
$
h(x) = {(10-15 e^{C}x+6 e^{2C}x^2)} / {(10-15x+6x^2)}
$
is decreasing in $x \in(0, e^{-C})$, for any constant $C\in(\log(4),
\infty)$.
\end{lemma}
\begin{pf}
The function $h(x)$ can be rewritten as
$
h(x) = \{6(-e^{C}x+5/4)^2+5/8\} / \{6(-x+5/4)^2+5/8\}$.
Note that $(-e^{C}x)/(-x)=e^{C}$ is nonincreasing in $x$
and $e^C>1$ for $x>0$.
Applying Lemma \ref{Lemma-2}, $(-e^{C}x+5/4)/(-x+5/4)$ is decreasing
in $x \in(0, e^{-C})$,
so is $(-e^{C}x+5/4)^2/(-x+5/4)^2$. When $C>\log(4)$,
${d\{(-e^{C}x+5/4)^2\}}/{dx} \le{d\{(-x+5/4)^2\}}/{dx}$.
This together with Lemma \ref{Lemma-2} verifies that
$h(x)$ is decreasing in $x\in(0,e^{-C})$.
\end{pf}

First, we evaluate $\alpha_{\infty}^{\FDR}$.
From (\ref{B.1}) and the conditions in Section \ref{sec-4.3},
%
%
\begin{equation}\label{C.1}
G_0(t)
= t \qquad\mbox{for } t\in[0,1]\quad \mbox{and}\quad
G_1(t)
=
\cases{
te^C, &\quad if $t\in[0, e^{-C}]$, \cr
1, &\quad if $t\in(e^{-C}, 1]$.}\hspace*{-28pt}
\end{equation}
Thus
$\sup_{0<t\le1} G_1(t)/t=e^C$. By $\widehat{\FDR}^{\infty}(t)$ in
Appendix \ref{appB},
%
%
\begin{equation} \label{C.2}
\alpha_{\infty}^{\FDR} =
\frac{\pi_0 +\pi_1 \{1-\lambda e^C\ID(\lambda<e^{-C})-\ID(\lambda
\ge e^{-C})\}/(1-\lambda)} {\pi_0+\pi_1 e^C}.
\end{equation}

Next, we compute $\alpha_{\infty}^{\FDR_L}$. Recall from Appendix \ref{appB} that
the distribution $G_0^*(t)$ with $k=5$ is that of a $\Betadist(3,3)$
random variable.
Similarly,\vspace*{1pt} by (\ref{C.1}),
the distribution $G_1^*(t)$ is that of a $\Betadist(3,3)/e^{C}$ random
variable.
By $\widehat\FDR_L^{\infty}(t)$ in Appendix~\ref{appB},
$\widehat\FDR_L^{\infty}(t)$ is a decreasing function of
$G_1^*(t)/G_0^*(t)$, for which two cases
need to be discussed.
In the first case, $t \in(0, e^{-C}]$, it follows that
\[
G_1^*(t)/G_0^*(t)
= e^{3C}\frac{10-15\cdot e^{C}t+6\cdot e^{2C}t^2}{10-15t+6t^2},
\]
which according to Lemma \ref{Lemma-3} is a decreasing function of $t$.
Thus,
\[
\sup_{t\in(0, e^{-C}]}G_1^*(t)/G_0^*(t) = \lim_{t\to0+}
G_1^*(t)/G_0^*(t) = e^{3C}
\]
and
%
%
\begin{equation} \label{C.3}
\inf_{t\in(0,e^{-C}]} \widehat{\FDR}_L^{\infty} (t)=
\frac{\pi_0+\pi_1\{1-G_1^*(\lambda)\}/\{1-G_0^*(\lambda)\}}
{\pi_0 +\pi_1 e^{3C}}.
\end{equation}
In the second case, $t\in(e^{-C},1]$, since $G_1^*(t) = 1$, we observe from
$\widehat\FDR_L^{\infty}(t)$ in Appendix \ref{appB} that
$\widehat\FDR_L^{\infty}(t)$
is an increasing function of $G_0^*(t)$, and thus
%
%
\begin{equation} \label{C.4}
\inf_{t\in(e^{-C}, 1]} \widehat{\FDR}_L^{\infty} (t)=
\frac{\pi_0+\pi_1\{1-G_1^*(\lambda)\}/\{1-G_0^*(\lambda)\}}
{\pi_0 +\pi_1 /G_0^*(e^{-C})}.
\end{equation}
Note that for $C>0$, we have
\[
\frac{1}{G_0^*(e^{-C})}=
\frac{e^{3C}}{6(e^{-C}-5/4)^2+5/8} \le
\frac{e^{3C}}{6(1-5/4)^2+5/8} = e^{3C}.
\]
Combining (\ref{C.3}) and (\ref{C.4}) gives
%
%
\begin{equation} \label{C.5}
\alpha_{\infty}^{\FDR_L} =
\frac{\pi_0+\pi_1\{1-G_1^*(\lambda)\}/\{1-G_0^*(\lambda)\}}
{\pi_0 +\pi_1 e^{3C}}.
\end{equation}
This completes the proof.
\end{appendix}

\section*{Acknowledgments}
The comments of the anonymous referees, the Associate Editor and the
Co-Editors are greatly appreciated.

\begin{supplement}
\stitle{Proofs and figures}
\slink[doi]{10.1214/10-AOS848SUPP}
\sdatatype{.pdf}
\sfilename{ZFY\_2\_supplement.pdf}
\sdescription{Section 1 gives detailed proofs of Theorems \ref{Thm-1}--\ref{Thm-3},
Section 2 gives the figure in Section \ref{sec-5.2}, and
Section 3 gives the figure in Section \ref{sec-5.3}.}
\end{supplement}

%
%


%
\printaddresses

\end{document}